\documentclass[]{article}
\usepackage{amsmath,amsfonts,amssymb}
\usepackage{amsthm}
\usepackage{mathtools}
\usepackage{cases}
\usepackage[mathfrak]{}
\usepackage[latin1]{inputenc}
\usepackage[T1]{fontenc}
\usepackage{verbatim}
\usepackage{graphicx}
\usepackage{float}
\usepackage{mathrsfs}
\usepackage [all]{xy}
\usepackage{color}
\usepackage{tikz}
\usepackage{comment}
\usepackage{hyperref}

\newtheorem{theorem}{Theorem}[section]
\newtheorem{lemma}[theorem]{Lemma}
\newtheorem{proposition}[theorem]{Proposition}%
\newtheorem{corollary}[theorem]{Corollary}

\newtheorem{example}{Example}%
\newtheorem{remark}{Remark}%

\newtheorem{definition}{Definition}

\usepackage{authblk}
\title{Sequences of point blow-ups over perfect fields from a combinatorial point of view}
\author[1]{Daniel Camaz\'on-Portela \footnote{The author was partially supported by grant PID2022-138906NB-C21 funded by MICIU/AEI/ 10.13039/501100011033 and by ERDF A way of making Europe.}}
\author[2]{Santiago Encinas \footnote{The author was partially supported by grant PID2022-138916NB-I00 funded by MCIN/AEI/10.13039/501100011033 and by ERDF A way of making Europe.}}
\affil[1]{Department of Mathematics, University of Almer\'ia, Carretera Sacramento, SN, Almer\'ia, 04120, Spain}
\affil[2]{Department of Algebra, Analysis, Geometry and Topology, University of Valladolid, Paseo Bel\'en 7, 47011, Spain}
\date{}                     
\begin{document}
  \maketitle



\begin{abstract}
We associate a combinatorial object to sequences of point blow-ups over perfect fields, the weighted directed graph, and another one to the composition of all blow-ups, which we call associated sequential morphisms, the $d-$ary intersection form. Then, in order to consider different fields extensions, we introduce the concepts of algebraically and combinatorially compatible partitions of the exceptional divisor for both sequences of point blow-ups and sequential morphisms, which lead us to define the corresponding algebraic and combinatorial equivalence classes. We prove that there exists a bijection between the respective combinatorial equivalence classes of sequences of point blow-ups and the associated sequential morphisms, and moreover, we also give a proof of the existence of a suitable bijection between the respective algebraic equivalence classes.
\end{abstract}

\section{Introduction}\label{sec1}
\setcounter{section}{4}

Sequences of blow-ups of smooth varieties along smooth centers are useful for general algebraic geometric purposes, in particular, for resolution and classification of singularities. We will assume, additionally, that the center of each blow-up has normal crossings with the already created exceptional divisors by the precedent blow ups, i.e. that at every point there are suitable regular systems of parameters such that both the center and the divisor components are locally defined by ideals generated by some of those parameters. 

Previously, sequences of blow-ups over perfect fields have been studied for different purposes. In \cite{BenitoVillamayor12}, Benito and Villamayor gave an overview of invariants of algebraic singularities over perfect fields and they showed how they lead to a synthetic proof of embedded resolution of singularities of $2-$dimensional schemes. Moreover, Benito and Villamayor in \cite{BenitoVillamayor13} continued the previous work in order to give invariants that would yield a sequence of monoidal transformations so as to eliminate the points of highest multiplicity of a hypersurface $X$. In \cite{ArvidssonBernasconiLacini22}, Arvidsson et al. proved the Kawamata-Viehweg vanishing theorem for surfaces of del Pezzo type over perfect fields of positive characteristic $p>5$, through the study of liftability of sequences of blow-ups of smooth points to characteristic $0$ over a smooth base. As a consequence, they showed that $klt$ threefold singularities over a perfect base field of characteristic $p>5$ are rational. Blache and Hallouin established a classification of singular del Pezzo surfaces over finite fields in \cite{BlacheHallouin23} which use sequences of blow-ups and contractions to describe most of the types. The same authors in \cite{BlacheHallouin24}, by taking advantage of \cite{BlacheHallouin23}, they selected eight types of non ordinary weak del Pezzo that are well suited for applications to error-correcting codes. In particular, they made use of blow-ups in order to add rational points and thus increase the length, as well as blow-downs of certain curves, those that are components of the most reducible sections of the anticanonical divisor, so as to improve the parameters of the code.

The composition of the successive blow ups of such sequences is a projective (and therefore proper) birational morphisms $Z\rightarrow X$, where $Z$ and $X$ are smooth algebraic varieties over a perfect field which are called sky and ground of the morphism, respectively. A morphism obtained in that way will be called sequential morphism. Each sequential morphism has an exceptional divisor which is the reduced divisor of $Z$ whose support consists of those points of $Z$ at which the morphism is not a local isomorphism. By the choice of the successive centers, the exceptional divisor of a sequential morphism has only normal crossing. In this paper, we restrict to the case of point blow-ups, as it is involved in many geometric contexts, in particular, the study of algebraic curves and surfaces. 
 
Campillo and Reguera in \cite{CampilloReguera94} studied sequences of point blow-ups using several combinatorial objects associated to the sequential morphism such as the d-ary intersection form on the abelian group of divisors with exceptional support, and also other equivalent objects as the weighted dual polyhedron or the weighted tree. These results are based on the proximity of infinitely near points. More details and other applications of proximity can be found on the survey \cite{CampilloGonzalezMonserrat09}. 


This paper is devoted to the study sequences of point blow-ups and the associated sequential morphisms defined over different algebraic extensions of a perfect field $k$. First, we associate a combinatorial object to sequential morphisms (the $d-$ary intersection form, see Definition \ref{DefdMultIntForm}) and another one to sequences of point blow-ups (the weighted directed graph, see Definition \ref{DefWeightdirectGraph}). Then, in order to consider different fields extensions, we introduce the key concepts of algebraically and combinatorially compatible partitions of the exceptional divisor for both sequences of point blow-ups and sequential morphisms. These concepts lead us to define the corresponding algebraic and combinatorial equivalence classes (see Definitions \ref{DefAlEqSM}, \ref{DefAlEqS}, \ref{ComEqSM} and \ref{Def3}). It follows that the combinatorial equivalence class of a weighted directed graph and that of the the $d-$ary intersection form are non-complete invariants for the class of algebraically equivalent sequences of point blow-ups and the class of algebraically equivalent sequential morphisms, respectively.

In order to establish a relation between these equivalence classes, we introduce the concept of final divisor of a sequential morphism: an irreducible component of the exceptional divisor is final, see Definition \ref{Def8}, if there is some other sequence of point blow-ups associated to the sequential morphism such that this irreducible component is the exceptional divisor of the last blow-up of this new sequence. Then, by using a geometrical approach, that makes an intensive use of intersection theory, we prove a result that is interesting by its own: a numerical characterization of final divisors.
\setcounter{theorem}{4}
\begin{proposition}\label{Pro7}
Let $(Z_{s},\ldots,Z_{0}, \pi)$ (resp. $\pi: Z_{s}\rightarrow Z_{0}$) be a sequence of point blow-ups (resp. sequential morphism) over a perfect field $K$ and let $H_{i}\in E_{i}^{s}$ be an irreducible component of the exceptional divisor of $\pi$. Then $H_{i}$ is final if and only if
\begin{equation*}
(H_{i})^{d}=(-1)^{r}(H_{i})^{s}\cdot (H_{j})^{r}\enspace \text{and}\enspace (H_{i})\cdot(H_{j})^{d-1}>0
\end{equation*}
for every $j$ such that $H_{i}\cap H_{j}\neq\emptyset$ (see Lemma \ref{Lem2} for a numerical characterization) and for all natural numbers $r$ and $s$ with $r+s=d$.
\end{proposition}
This characterization allow us to show that the intersection form on the sky $Z=Z_{s}$ determines those exceptional components which are final. By making use of the projection formula, the intersection form on the blow-down of the final components can be explicitly computed in terms of the intersection form of $Z$. By focusing on birational contractions (see \cite{Ishii77}), we prove that a regular projective contraction of a final divisor to a point verifies that the image of the exceptional divisors is still a simple normal crossing divisor and that the variety obtained corresponds to the sky of a sequence of point blow-ups with ground $Z_{0}$.
\setcounter{theorem}{10}
\begin{proposition}\label{Pro8}
Let $(Z_{s},\ldots,Z_{0}, \pi)$ be a sequence of point blow-ups (as in Definition \ref{Def1}) of length $m$ and let $H_{i}\in E_{i}^{s}$ be an irreducible component of the exceptional divisor of $\pi$. If $H_{i}$ is final, then there exists a regular projective contraction $(Z, f_{m}, X_{m-1})$ of $H_{i}$ to a point such that $f_{m}(E_{Z})$ is a simple normal crossing divisor and $X_{m-1}$ is the sky of a sequence of point blow-ups with ground $Z_{0}$.
\end{proposition}
This lead us to establish a well defined bijection between combinatorial equivalence classes of weighted directed graphs of sequences of point blow-ups and $d-$ary intersection forms of the associated sequential morphism.
\setcounter{theorem}{11}
\begin{theorem}\label{Thm1}
Let $(Z_{s},...,Z_{0},\sqcup_{i=1}^{l}F_{i},\pi)_{comb}$ and $(Z^{'}_{s},...,Z^{'}_{0},\sqcup_{i=1}^{l^{'}}F_{i}^{'},\pi^{'})_{comb}$ be two combinatorially marked sequences of point blow-ups, with $l=l^{'}$. The associated weighted directed graphs $T_{(Z_{s},...,Z_{0}, \sqcup_{i=1}^{l}F_{i}, \pi)}$, $T_{(Z_{s}^{'},...,Z_{0}^{'}, \sqcup_{i=1}^{l^{'}}F_{i}^{'}, \pi^{'})}$ are equivalent if and only if the multilinear maps $\Phi_{Z, \left(\sqcup_{i=1}^{l}F_{i}\right)_{Z, \pi}, K}$ and $\Phi_{Z^{'}, \left(\sqcup_{i=1}^{l}F_{i}^{'}\right)_{Z^{'}, \pi^{'}}, K}$ associated to the corresponding combinatorially marked sequential morphisms are also equivalent as in Definition \ref{ComEqSM}.

\end{theorem}
Moreover, we go an step further and we prove that is is possible to recover from the sequential morphism, up to algebraic equivalence, the sequence of point blow-ups itself.
\setcounter{theorem}{13}
\begin{proposition}\label{Thm41}
Let $(\pi: Z_{s}\rightarrow Z_{0},\sqcup_{i=1}^{l}F_{i})_{alg}$ be an algebraically marked sequential morphism. Given the $d-$ary multilinear intersection form associated to the partition, $\mathcal{I}_{Z, \left(\sqcup_{i=1}^{l}F_{i}\right)_{Z, \pi}, K}$ (see Definition \ref{MInFPar}), we can recover all the algebraically marked sequences of point blow-ups that are associated to algebraically marked sequential morphisms in the same algebraic equivalence class of $(\pi: Z_{s}\rightarrow Z_{0},\sqcup_{i=1}^{l}F_{i})_{alg}$.
\end{proposition}
This result allow us to establish a bijection between algebraic equivalence classes of sequences of point blow-ups and sequential morphisms. 



\setcounter{theorem}{16}
\begin{theorem}\label{Thm21}
Let $(\pi: Z_{s}\rightarrow Z_{0}, \sqcup_{i=1}^{l}F_{i})_{alg}$ and $(\pi^{'}: Z_{s}^{'}\rightarrow Z_{0}^{'}, \sqcup_{i=1}^{l^{'}}F^{'}_{i})_{alg}$ be two algebraically marked sequential morphisms. Then  they are algebraically equivalent over $K$ as in Definition \ref{DefAlEqSM} if and only if
 there exist algebraically marked sequences of point blow-ups $(Z_{s},...,Z_{0},\sqcup_{i=1}^{l}F_{i},\pi)_{alg}$ and $(Z^{'}_{s},...,Z^{'}_{0},\sqcup_{i=1}^{l^{'}}F_{i}^{'},\pi^{'})_{alg}$  associated to $(\pi: Z_{s}\rightarrow Z_{0}, \sqcup_{i=1}^{l}F_{i})_{alg}$ and $(\pi^{'}: Z_{s}^{'}\rightarrow Z_{0}^{'}, \sqcup_{i=1}^{l^{'}}F^{'}_{i})_{alg}$ respectively such that they are algebraically equivalent over $K$ as in Definition \ref{DefAlEqS}.
\end{theorem}

Finally, we illustrate within Example \ref{Exmp1} that the numerical information contained in the $d-$ary intersection form is minimal in some sense for our classification purpose.



\setcounter{section}{1}
\section{Preliminaries}

\subsection{Sequences of blow-ups and sequential morphisms}

\setcounter{theorem}{0}

Fix a perfect field $k$ and chose an algebraic closure $\overline{k}$. Throughout this paper, unless otherwise stated, a variety will mean a reduced projective scheme over a perfect field $K$, with $K$ and algebraic extension of $k$ such that $K\subset\overline{k}$, and a point will mean a closed point.

\begin{definition}\label{Def1}
A sequence of blow-ups over $K$ is defined as a sequence of blow-ups at smooth closed subvarieties $C_{i}$ of smooth $d-$dimensional projective varieties $Z_{i}$
\begin{equation*}
Z=Z_{s}\xrightarrow{\pi_{s}} Z_{s-1}\xrightarrow{\pi_{s-1}}\cdot\cdot\cdot\xrightarrow{\pi_{2}} Z_{1}\xrightarrow{\pi_{1}} Z_{0}=X,
\end{equation*}
such that for $i\in\left\{0,1,...,s-1\right\}$:
\begin{enumerate}
\item if we denote by $C_{i+1}$ to the center of $\pi_{i+1}$, then $C_{i+1}$ is a smooth subvariety of $Z_{i}$ defined over $K$ that could be reducible, that is $C_{i+1}=\sqcup C_{i+1,j}$ with $C_{i+1,j}$ irreducible over $K$
\item $codim(C_{i+1})\geq 2$
\item if we denote by $E_{j}^{j}$ the exceptional divisor of $\pi_{j}$, and for $k>j$ we denote by $E_{j}^{k}$ the strict transform of $E_{j}^{j}$ in $Z_{k}$, then $C_{i+1}$ has simple normal crossings with $\{E_{1}^{i}, E_{2}^{i},...,E_{i}^{i}\}$
\end{enumerate}
\end{definition}

We denote by $\pi$ the composition $\pi_{1}\circ\pi_{2}\circ...\circ\pi_{s-1}\circ\pi_{s}$.

\begin{definition}\label{DefSM}
A morphism $\pi: Z\rightarrow Z_{0}$ which can be expressed, in at least one way, as a composition of blow ups with the conditions in Definition \ref{Def1} will be called a sequential morphism over $K$.
\end{definition}

\begin{remark}\label{NoPi}
Given a sequence of blow-ups $\pi$, we denote by $\pi_{s,i}: Z\rightarrow Z_{i}$ where $\pi_{s,i}=\pi_{i+1}\circ\pi_{i+2}\circ...\circ\pi_{s-1}\circ\pi_{s}$. 
\end{remark}

\begin{remark}
We will refer to $X=Z_{0}$ and $Z$ as the ground and the sky of the sequential morphism $\pi$ respectively.
\end{remark}

\begin{definition}
The length $m$ over $K$ of a sequence of blow ups is defined as $\sum_{i=1}^{s}\# C_{i}$, where $\# C_{i}$ denotes the number of irreducible components of $C_{i}$ over $K$. Notice that it coincides with the number of irreducible components of the exceptional divisor $E$ (over $K$ too). Therefore, the length depends on the sequential morphism $\pi$ and it can be also called the length of $\pi$ over $K$. It will be denoted by $m=lenght_{K}(\pi)$. Notice that $s\leq m$, and $s=m$ exactly when all the blow up centers are irreducible over $K$.
\end{definition}

\begin{remark}
Note that in the case of sequences of point blow-ups, that is $C_{i}=P_{i}$, if $K=\overline{k}$, with $\overline{k}$ the algebraic closure of $k$, $m=lenght_{K}(\pi)=\sum_{i=1}^{s}\left[K(P_{i}):K\right]$, where $K(P_{i})$ is the residue field of $P_{i}$.
\end{remark}

\begin{remark}
Moreover we will denote by $H_{\beta}$ the irreducible components over $K$ of the exceptional divisor $E$ of $\pi$, that is we have $E=\sqcup_{\beta} H_{\beta}$.
\end{remark}

\begin{definition}\label{DefProx}
Let $(Z_{s},...,Z_{0}, \pi)$ be a sequence of blow-ups of length $m$, and $H_{1},...,H_{m}$ the irreducible components of the exceptional divisor $E$ over $K$ of the associated sequential morphism. For each $i$, with $i=1,2,..,m$, let $r(i)$ be the integer such that the image of $H_{i}$ at $Z_{r(i)}$ is a component of the center (codimension at least $2$) whose blow-up creates $H_{i}$. If $j$ is different from $i$, and the image of $H_{j}$ at $Z_{r(i)}$ has codimension $1$ and it contains the image of $H_{i}$ at $Z_{r(i)}$, then $H_{i}$ is said to be proximate to $j$, and we denote it by $H_{i}\rightarrow H_{j}$. It is clear that one has $r(i)>r(j)$ when $H_{i}$ is proximate to $H_{j}$.
\end{definition}

For sequences of point blow-ups we denote $deg(H_{i})=\left[K(P_{i}):K\right]$, where $P_{i}$ is the point in the center of $\pi_{r(i)}$ such that the image of $H_{i}$ in $Z_{r(i)}$ is $P_{i}$.

\subsection{Some useful intersection theory}

In this section we will recall some well known results on intersection theory which we will use in the proof of our results. We do not presents an exhaustive list, but the aim of this section is to list some results which either are critical in our development, or are very specific in intersection theory, or we did not find any proof in the literature.

Along this section, $X$ and $Z$ will stand for smooth algebraic varieties over $K$, the considered maps $Z\rightarrow X$ for proper morphisms, and $d$ for the dimension of $X$.

\begin{remark}
Let $V$ be a $n-$dimensional irreducible subvariety of $X$. Although it is common to denote by $\left[V\right]$ to the equivalence class of $V$ in the Chow ring $A^{\bullet}(X)$, for simplicity we will denote also by $V$ to the equivalence class whenever there is not possible confusion.
\end{remark}

\begin{definition}\label{Def6}\cite[Section 2.6]{Fulton98}
Let $D$ be an effective Cartier divisor of $X$, and let $i: D\rightarrow X$ be the inclusion. There are Gysin homorphisms:
\begin{equation*}
i^{*}: Z_{j}X\rightarrow A_{j-1}D,
\end{equation*}
for $j=1, \ldots, dim(X)$, defined by the formula
\begin{equation*}
i^{*}(\alpha)=D\cdot\alpha.
\end{equation*}
\end{definition}

\begin{proposition}\label{Pro1}\cite[Proposition 2.6.]{Fulton98}
There are therefore induced homorphisms:
\begin{equation*}
i^{*}: A_{j}X\rightarrow A_{j-1}D
\end{equation*}
for $j=1, \ldots, dim(X)$ such that one has
\renewcommand{\labelenumi}{\alph{enumi}}
\begin{enumerate}
\item If $\alpha$ is a $j$-cycle on $X$, then
\begin{equation*}
i_{*}i^{*}(\alpha)=c_{1}(\mathcal{O}_{X}(D))\cap\alpha;
\end{equation*}
\item If $\alpha$ is a $j$-cycle on $D$, then
\begin{equation*}
i^{*}i_{*}(\alpha)=c_{1}(N_{D/X})\cap\alpha.
\end{equation*}
\end{enumerate}
\end{proposition}

Notice that if $f:Z\rightarrow X$ is a proper morphism, $V$ is an irreducible subvariety of $Z$ and $W=f(V)$ its image by $f$, then one defines the push-forward of the cycle associated to $V$ as follows:
\begin{equation}\label{EqPushForw}
f_{*}(V)=
\begin{cases}
0 & \text{if}\enspace dim V>dim W, \\
[K(V):K(W)]W & \text{if}\enspace dim V=dim W, 
\end{cases}
\end{equation}
where $[K(V):K(W)]$ denotes the degree of the extension of function fields $K(V)/K(W)$.\\
This definition allows also to define, by linearity, the push-forward $f_{*}(\alpha)$ for any cycle $\alpha$.

The next result, known as the projection formula, is specially useful and it will play an essential role in some of our intersection computations.
\begin{proposition}\label{Pro2}\cite[Proposition2.3.]{Fulton98}
Let $D$ be a divisor on $X$, $f: Z\rightarrow X$ a proper morphism, $\alpha$ a $k-$cycle on $Z$, and $g$ the morphism from $f^{-1}(\left|D\right|)\cap\left|\alpha\right|$ to $\left|D\right|\cap f(\left|\alpha\right|)$ induced by $f$. Then
\begin{equation*}
g_{*}(f^{*}D\cdot\alpha)=D\cdot f_{*}(\alpha)
\end{equation*}
in $A_{k-1}(\left|D\right|\cap f(\left|\alpha\right|))$.
\end{proposition}

The following proposition is obtained as a particular case of the one above.
\begin{proposition}\label{Pro3}\cite[Proposition 1.10]{Debarre01}
Let $f: Z\rightarrow X$ be a proper surjective morphism. Let $D_{1}, D_{2},..., D_{r}$ be Cartier divisors on $X$ with $r=d=dim(X)$. Then, one has
\begin{equation*}
f^{*}D_{1}\cdot f^{*}D_{2}\cdot\cdot\cdot f^{*}D_{r}=deg(f)D_{1}\cdot D_{2}\cdot\cdot\cdot D_{r}
\end{equation*}
where $deg(f)=\left[K(Z):K(X)\right]$, if $deg(f)$ is finite.
\end{proposition}

Now, let $P$ be a closed point of a variety $X$, $\widetilde{X}$ the blow-up of $X$ at $P$, and $E$ the exceptional divisor. Then, we have the following fibre square:
\begin{equation*}
\xymatrix{
E\ar[r]^{j}\ar[d]^{g} & \widetilde{X}\ar[d]_{f}\\
P\ar[r]^{i} & X
}
\end{equation*}
In order to apply the projection formula (either Proposition \ref{Pro2} or \ref{Pro3}), then we need previously the following result.
\begin{theorem}\cite[Corollary 6.7.1.]{Fulton98}\label{ThmTotTrans}
Let $V$ be a $r-$dimensional subvariety of $X$, and let $\widetilde{V}\subset\widetilde{X}$ be the proper transform of $V$. Then
\begin{equation*}
f^{*}(V)=\widetilde{V}+(e_{P}V)j_{*}(L),
\end{equation*}
in $A_{r}(\widetilde{X})$, where $L$ is a $r-$dimensional linear subspace of $E=\mathbb{P}_{K(P)}^{d-1}$, $K(P)$ is the residue field of $\mathcal{O}_{P,X}$, and $e_{P}V$ is the multiplicity of $P$ on $V$.
\end{theorem}

As it is shown in Proposition \ref{Pro1}, the normal bundle plays a key role in many intersection computations. The next two results refers to the computation of the normal bundle of certain varieties we are interested in. More concretely, varieties resulting from the intersection of two Cartier divisors meeting regularly (i.e. transversally at smooth points of both) and varieties obtained as the strict transforms of other ones by a blow-up.
\begin{lemma}\label{Lem1}
Let $D$ and $F$ be two irreducible components of a simple normal crossing divisor $E$ that is regularly embedded in $X$. If we denote by $G=D\cap F$ then
\begin{equation*}
N_{F/X}\vert_{G}\cong N_{G/D}.
\end{equation*}
\end{lemma}

\begin{proof}
Let $i_{G, D}: G\rightarrow D$ and $i_{D,X}: D\rightarrow X$ be regular embeddings. Then the composite $i_{D, X}\circ i_{G, D}$ is a regular imbedding, and there is an exact sequence of vector bundles on $G$ (see \cite[Proposition 19.1.5]{EGAIV})
\begin{equation*}
0\rightarrow N_{G/D}\rightarrow N_{G/X}\rightarrow N_{D/X}\vert_{G}\rightarrow 0.
\end{equation*}  
Since $D$ and $F$ meet regularly in $X$, then (see \cite[Proposition 3.6.]{Fulton85})
\begin{equation*}
N_{G/X}\cong N_{D/X}\vert_{G}\oplus N_{F/X}\vert_{G},
\end{equation*}
so we have the following exact sequences of vector bundles on $G$
\begin{equation*}
0\rightarrow N_{G/D}\rightarrow N_{D/X}\vert_{G}\oplus N_{F/X}\vert_{G}\rightarrow N_{D/X}\vert_{G}\rightarrow 0
\end{equation*}
Then it follows that $N_{G/D}\cong N_{F/X}\vert_{G}$.
\end{proof}

\begin{proposition}\label{Pro4}\cite[Appendix B.6.10.]{Fulton98}
If $C\subset Y$ and $Y\subset X$ are regular embeddings, let $\widetilde{X}=Bl_{C}X$, $E$ the exceptional divisor in $\widetilde{X}$, $\rho$ the projection from $\widetilde{X}$ to $X$. Let $\widetilde{Y}=Bl_{C}Y$. Then $\widetilde{Y}\subset\rho^{-1}(Y)$, $E\subset\rho^{-1}(Y)$, and $\widetilde{Y}$ is the residual scheme to $E$ in $\rho^{-1}(Y)$, i.e., the ideal sheaves of $\widetilde{Y}$, $E$ and $\rho^{-1}(Y)$ in $\widetilde{X}$ are related by
\begin{equation*}
\mathcal{I}(\widetilde{Y})\cdot\mathcal{I}(E)=\mathcal{I}(\rho^{-1}(Y)).
\end{equation*}
In addition, the canonical embedding of $\widetilde{Y}$ in $\widetilde{X}$ is a regular imbedding, with normal bundle
\begin{equation*}
N_{\widetilde{Y}/\widetilde{X}}\cong\pi^{*}N_{Y/X}\otimes\mathcal{O}(-F)
\end{equation*}
where $\pi$ is the projection from $\widetilde{Y}$ to $Y$, and $F$ is the exceptional divisor on $\widetilde{Y}$ of such projection.
\end{proposition}

\begin{proposition}\label{Pro5}\cite[Appendix B.6.3]{Fulton98}
With the same notations as above
\begin{equation*}
N_{E/\widetilde{Z}}\cong\mathcal{O}_{E}(-1).
\end{equation*}
\end{proposition}

Finally, last results of this section deals with the structure of the Chow ring of a projective bundle over a smooth projective variety and the computation of the class corresponding to a projective subbundle. These are crucial for the rest of the paper as they describe the intersection structure of the irreducible components of the exceptional divisor obtained after a blow-up.
\begin{theorem}\label{Thm3}\cite[Theorem 9.6.]{EisenbudHarris16}
Let $\mathcal{E}$ be a vector bundle of rank $r+1$ on a smooth projective variety $X$, and let $\varsigma=c_{1}(\mathcal{O}_{\mathbb{P}\mathcal{E}}(1))\in A^{1}(\mathbb{P}\mathcal{E})$, and $p_{\mathbb{P}\mathcal{E}}: \mathbb{P}\mathcal{E}\rightarrow X$ the projection of the induced projective bundle. The map $p_{\mathbb{P}\mathcal{E}}^{*}: A(X)\rightarrow A(\mathbb{P}\mathcal{E})$ is an injective ring homomorphism , and via this map one has the isomorphism of A(X)-algebras given by
\begin{equation*}
A(\mathbb{P}\mathcal{E})\cong A(X)\left[\varsigma\right]/\varsigma^{r+1}+c_{1}(\mathcal{E})\varsigma^{r}+\cdot\cdot\cdot+c_{r+1}(\mathcal{E})
\end{equation*}
\end{theorem}

\begin{corollary}\cite[Theorem 2.1.]{EisenbudHarris16}
The Chow ring of $\mathbb{P}^{n}$ is
\begin{equation*}
A(\mathbb{P}^{n})\cong\mathbb{Z}\left[\varsigma\right]/(\varsigma^{n+1})
\end{equation*}
where $\varsigma\in A^{1}(\mathbb{P}^{n})$ is the rational equivalence class of all the hyperplane; more generally, the
class of all closed subschemes of codimension $k$ and degree $d$ is $d\varsigma^{k}$
\end{corollary}

\begin{proposition}\label{Pro6}\cite[Proposition 9.13]{EisenbudHarris16}
If $X$ is a smooth projective variety and $\mathcal{F}\subset\mathcal{E}$ are vector bundles on $X$ of ranks $s$ and $r$ respectively, then
\begin{equation*}
\left[\mathbb{P}(\mathcal{F})\right]=\varsigma^{r-s}+\gamma_{1}\varsigma^{r-s-1}+...+\gamma_{r-s}\in A^{r-s}(\mathbb{P}(\mathcal{E}))
\end{equation*}
, where $\varsigma=c_{1}(\mathcal{O}_{\mathbb{P}(\mathcal{E})}(1))$ and $\gamma_{k}=c_{k}(\mathcal{E}/\mathcal{F})$. Moreover, the normal bundle of $\mathbb{P}\mathcal{F}$ in $\mathbb{P}\mathcal{E}$ is $\mathcal{O}_{\mathbb{P}\mathcal{F}}(1)\otimes p_{\mathbb{P}\mathcal{F}}^{*}(\mathcal{E}/\mathcal{F})$.
\end{proposition}

\subsection{Combinatorial objects associated to sequences of point blow-ups and sequential morphisms}

This section is devoted to define two combinatorial objects associated to sequences of point blow-ups and sequential morsphisms: the weighted directed tree and the $d-$ary intersection form, respectively.
\begin{definition}
We denote by $E_{Z, \pi, K}$ (or simply $E_{Z}$ if there is no confusion) the abelian group of divisors of $Z$ of the form $\sum_{\beta=1}^{m}n_{\beta}H_{\beta}$ where $n_{\beta}\in\mathbb{Z}$. In other words,$E_{Z, \pi, K}$ is the free group generated by the irreducible components of $E$ over $K$. The set $\left\{H_{\beta}\right\}_{\beta=1}^{m}$ is a basis of the $\mathbb{Z}-$module $E_{Z}$.
\end{definition}

\begin{definition}\label{DefdMultIntForm}
Given a sequential morphism $\pi: Z\rightarrow Z_{0}$, we consider the $d$-ary multilinear intersection form:
\begin{equation*}
\mathcal{I}_{Z, E_{Z, \pi}, K}: \overbrace{E_{Z}\times E_{Z}\times\cdots\times E_{Z}}^{d}\rightarrow \mathbb{Z},
\end{equation*}
defined by intersecting cycles in the sky $Z$ and taking degrees, that is,
\begin{equation*}
\mathcal{I}_{Z, E_{Z, \pi}, K}(H_{i_{1}}, H_{i_{2}},..., H_{i_{d}})=deg(H_{i_{1}}\cdot H_{i_{2}}\cdot H_{i_{3}}\cdots H_{i_{d}}),
\end{equation*}
where $H_{i_{1}}\cdot H_{i_{2}}\cdot H_{i_{3}}\cdots H_{i_{d}}$ is an intersection class of $0-$cycles in the abelian group $A_{0}(Z)$, and $deg$ stands for the degree in the Chow ring $A(Z)$.
\end{definition}

\begin{remark}
For the sake of simplicity we will denote $H_{i_{1}}\cdot H_{i_{2}}\cdots H_{i_{d}}=\mathcal{I}_{Z, E_{Z, \pi}, K}(H_{i_{1}}, H_{i_{2}},..., H_{i_{d}})$.
\end{remark}

\begin{definition}
Given a sequential morphism $\pi$ as in Definition \ref{DefSM}, it induces a natural isomorphism $E_{Z}\cong\mathbb{Z}^{m}$, where the standard basis of $\mathbb{Z}^{m}$ is the image of the $\mathbb{Z}$-basis $\left\{H_{i}\right\}_{i=1}^{m}$. In this way, the $d-$ary interection form gives rise to a multilinear map:
\begin{equation*}
\Phi_{Z, E_{Z, \pi}, K}:\overbrace{\mathbb{Z}^{m}\times\cdot\cdot\cdot\times\mathbb{Z}^{m}}^{d}\rightarrow\mathbb{Z}.
\end{equation*}
We say that $\Phi_{Z, E_{Z, \pi}, K}$ is the multilinear map associated to $\pi$. The permutation group $\mathcal{S}_{m}$ acts on the set of multilinear forms $\mathbb{Z}^{m}\rightarrow\mathbb{Z}$  by permuting the elements of the standard basis of $\mathbb{Z}^{m}$. It is clear that if we denote by $\Psi_{Z, E_{Z, \pi}, K}$ to the orbit of $\Phi_{Z, E_{Z, \pi}, K}$, then $\Psi_{Z, E_{Z, \pi}, K}$ does not depend on the labeling of the elements of the basis $\left\{H_{i}\right\}$.
\end{definition}

\begin{definition}\label{DefWeightdirectGraph}
Given a sequence of point blow-ups $(Z_{s},\ldots,Z_{0}, \pi)$ as in Definition \ref{Def1}, we associate to it a directed graph $T_{(Z_{s},\ldots,Z_{0}, \pi)}$ in the following way: The vertices of $T_{(Z_{s},\ldots,Z_{0}, \pi)}$ are in one to one correspondence with the points $P_{1}, P_{2}, \ldots, P_{s}$ corresponding to the centers, and the edges with the pairs $(P_{i}, P_{j})$ such that $j\geq i+1$ and $P_{j}$ is proximate to $P_{i}$. For each $P_{i}$, we denote by $\mathit{l}_{i}$ the degree of the finite field extension $K\subset K(P_{i})$, where $K(P_{i})$ is the residue field of $P_{i}$. We define the weighted directed graph associated to the sequence $(Z_{s},\ldots,Z_{0}, \pi)$  as the directed graph together with the following weights on its vertices: The vertex corresponding to the point $P_{i}$ has the value $\mathit{l}_{i}$ as weight.
\end{definition}


\section{Compatible partitions of the exceptional divisor and its associated equivalence classes}

In order to consider different fields $K$, with $k\subset K\subset\overline{k}$, we define the notion of a compatible partition of the exceptional divisor $E$. This concept will be set for sequential morphisms (Definitions \ref{AlComPartSM} and \ref{CombCompParSM}) and sequences of blow-ups (Definitions \ref{AlComPartSeq} and \ref{CombCompParSeq}.

\begin{definition}\label{AlComPartSM}
Given a sequential morphism $\pi: Z\rightarrow Z_{0}$ as in Definition \ref{Def1}, and a partition of the exceptional divisor $E$, $\sqcup_{i=1}^{l}F_{i}$, we will say that the partition is algebraically compatible with the morphism $\pi$ if there exists a smaller field $\widetilde{K}\subset K$ with $k\subset\widetilde{K}$, there are $\widetilde{K}-$varieties $\widetilde{Z_{0}}$ and $\widetilde{Z}$ and a $\widetilde{K}-$morphism $\widetilde{Z}\xrightarrow{\widetilde{\pi}}\widetilde{Z_{0}}$:
\begin{equation*}
\xymatrix{
Z\cong\widetilde{Z}\times_{Spec(\widetilde{K})} Spec(K)\ar[r]^{\pi}\ar[d]_{\beta} & Z_{0}\cong\widetilde{Z_{0}}\times_{Spec(\widetilde{K})} Spec(K)\ar[d] \\
\widetilde{Z} \ar[r]^{\widetilde{\pi}} & \widetilde{Z_{0}},
}
\end{equation*}
such that if we denote by $\widetilde{H_{1}}, \ldots,\widetilde{H_{l}}$ the irreducible components of the exceptional divisor $\widetilde{E}$, then, for each $i=1,...,l$, $\beta(H)=\widetilde{H_{i}}$ for all $H\in F_{i}$.
\end{definition}
An algebraically marked sequential morphism is defined by a pair $(\pi: Z\rightarrow Z_{0},\sqcup_{i=1}^{l}F_{i})_{alg}$, where $ \sqcup_{i=1}^{l}F_{i}$ is a partition algebraically compatible with the sequential morphism $\pi$. This notion leads in a natural way to the following equivalence relation.
\begin{definition}\label{DefAlEqSM}
We say that two algebraically marked sequential morphisms $(\pi: Z\rightarrow Z_{0},\sqcup_{i=1}^{l}F_{i})_{alg}$ and $(\pi^{'}: Z^{'}\rightarrow Z_{0}^{'}, \sqcup_{i=1}^{l^{'}}F^{'}_{i})_{alg}$ over $K$ are algebraically equivalent, and we denote it by $(\pi: Z\rightarrow Z_{0},\sqcup_{i=1}^{l}F_{i})_{alg}\sim(\pi^{'}: Z^{'}\rightarrow Z_{0}^{'}, \sqcup_{i=1}^{l^{'}}F^{'}_{i})_{alg}$ , if and only if there exist smaller fields $\widetilde{K},\widetilde{K^{'}}\subset K$ with $\widetilde{K}\cong_{k}\widetilde{K^{'}}$ satisfying the conditions of definition \ref{AlComPartSM}:
\begin{equation*}
\xymatrix{
Z\cong\widetilde{Z}\times_{Spec(\widetilde{K})} Spec(K)\ar[r]^{\pi}\ar[d] & Z_{0}\cong\widetilde{Z_{0}}\times_{Spec(\widetilde{K})} Spec(K)\ar[d] \\
\widetilde{Z} \ar[r]^{\widetilde{\pi}} & \widetilde{Z_{0}}
}
\end{equation*}
\begin{equation*}
\xymatrix{
Z^{'}\cong\widetilde{Z^{'}}\times_{Spec(\widetilde{K^{'}})} Spec(K)\ar[r]^{\pi^{'}}\ar[d] & Z^{'}_{0}\cong\widetilde{Z^{'}_{0}}\times_{Spec(\widetilde{K^{'}})} Spec(K)\ar[d] \\
\widetilde{Z^{'}} \ar[r]^{\widetilde{\pi^{'}}} & \widetilde{Z^{'}_{0}}
}
\end{equation*}
and there exist isomorphisms $a: \widetilde{Z_{0}}\rightarrow\widetilde{Z_{0}^{'}}$ and $b: \widetilde{Z}\rightarrow\widetilde{Z^{'}}$ such that the following diagram is commutative:
\begin{equation*}
\xymatrix{\widetilde{Z}\ar[r]^{b}\ar[d]^{\widetilde{\pi}} & \widetilde{Z^{'}}\ar[l]\ar[d]^{\widetilde{\pi^{'}}}\\
\widetilde{Z_{0}}\ar[r]^{a} & \widetilde{Z_{0}^{'}}\ar[l]}
\end{equation*}
\end{definition}


\begin{definition}\label{AlComPartSeq}
Given a sequence of point blow-ups $(Z_{0},...,Z_{s},\pi)$ and a partition of the exceptional divisor $E$, $\sqcup_{i=1}^{l}F_{i}$, we sill say that the partition is algebraically compatible with the sequence $(Z_{0},...,Z_{s},\pi)$ if there exists a smaller field $\widetilde{K}\subset K$ with $k\subset\widetilde{K}$ and there are $\widetilde{K}-$varieties $\widetilde{Z_{i}}$ and $\widetilde{K}-$morphisms $\widetilde{Z}_{i+1}\xrightarrow{\widetilde{\pi}_{i+1}}\widetilde{Z}_{i}$
\begin{equation*}
\xymatrix{
Z_{s}\ar[r]\ar[d] & Z_{s-1}\ar[r]\ar[d] & \cdot\ar[r] & Z_{1}\ar[r]\ar[d] & Z_{0}\ar[d] \\
\widetilde{Z}_{s}\ar[r] & \widetilde{Z}_{s-1}\ar[r] & \cdot\ar[r] & \widetilde{Z}_{1}\ar[r] & \widetilde{Z}_{0}
}
\end{equation*}
where $Z_{i}\cong\widetilde{Z}_{i}\times_{Spec(\widetilde{K})} Spec(K)$ for all $i=1, \ldots,s$ , such that the exceptional divisor of $(\widetilde{Z_{0}},...,\widetilde{Z_{l}},\widetilde{\pi})$ has irreducible components $\widetilde{H_{1}},...,\widetilde{H_{l}}$ and, for each $i=1,...,l$, $\beta(H)=\tilde{H}_i$ for all $H\in F_i$.
\end{definition}

\begin{remark}
Notice that, since $k$ is perfect, $\tilde{K}\subseteq K$ is a separable algebraic extension.
\end{remark}

An algebraically marked sequence of point blow-ups is defined by a pair $(Z_{s},...,Z_{0},\pi,\sqcup_{i=1}^{l}F_{i})_{alg}$, where $ \sqcup_{i=1}^{l}F_{i}$ is a partition algebraically compatible with $(Z_{s},...,Z_{0},\pi)$. Before defining the corresponding equivalence relation, we need to introduce the concept of brick blow-up.
\begin{definition}\label{DefBriBl}
Given a  variety $X$ we will call a brick blow-up with ground $X$ to a sequential morphism obtained as a composition of point blow-ups with disjoint centers $\sqcup_{j=1}^{l} C_{j}\subset X$, $X^{'}=X_{l}\rightarrow X_{l-1}\rightarrow ...\rightarrow X_{1}\rightarrow X$. Note that $X_{i}\rightarrow X_{i-1}$ is the brick blow-up at $C_{i}$, where $C_{i}$ need not to be irreducible.
\end{definition}

\begin{definition}\label{DefAlEqS}
We say that two algebraically marked sequences of point blow ups, $(Z_{s},...,Z_{0},\pi,\sqcup_{i=1}^{l}F_{i})_{alg}$ , and $(Z_{s^{'}},...,Z^{'}_{0},\pi^{'},\sqcup_{i=1}^{l^{'}}F_{i}^{'})_{alg}$, are algebraically equivalent over $K$, and we denote it by $(Z_{s},...,Z_{0},\pi,\sqcup_{i=1}^{l}F_{i})_{alg}\sim(Z_{s^{'}},...,Z^{'}_{0},\pi^{'},\sqcup_{i=1}^{l^{'}}F_{i}^{'})_{alg}$, if $m=m^{'}$ and there exist smaller fields $\widetilde{K},\widetilde{K^{'}}\subset K$ with $\widetilde{K}\cong_{k}\widetilde{K^{'}}$ 
\begin{equation*}
\xymatrix{
Z_{s}\ar[r]\ar[d] & Z_{s-1}\ar[r]\ar[d] & \cdot\ar[r]\ar[d] & Z_{1}\ar[r]\ar[d] & Z_{0}\ar[d] \\
\widetilde{Z}_{s}\ar[r] & \widetilde{Z}_{s-1}\ar[r] & \cdot\ar[r] & \widetilde{Z}_{1}\ar[r] & \widetilde{Z}_{0}
}
\end{equation*}
\begin{equation*}
\xymatrix{
Z^{'}_{s}\ar[r]\ar[d] & Z^{'}_{s-1}\ar[r]\ar[d] & \cdot\ar[r]\ar[d] & Z^{'}_{1}\ar[r]\ar[d] & Z^{'}_{0}\ar[d] \\
\widetilde{Z}^{'}_{s}\ar[r] & \widetilde{Z}^{'}_{s-1}\ar[r] & \cdot\ar[r] & \widetilde{Z}^{'}_{1}\ar[r] & \widetilde{Z}^{'}_{0}
}
\end{equation*}
 with $Z_{i}\cong\widetilde{Z}_{i}\times_{Spec(\widetilde{K})} Spec(K)$ (resp. $Z_{i}^{'}\cong\widetilde{Z}_{i}^{'}\times_{Spec(\widetilde{K}^{'})} Spec(K)$) and algebraic isomorphisms $a, b=b_{t}, b_{t-1},...,b_{1}$, with $t\leq s$, such that there are indexes $r_{1},...,r_{t}=s\in\left\{1,...,l\right\}$ and $r_{1}^{'},...,r_{t}^{'}=s^{'}\in\left\{1,...,s^{'}\right\}$, where $Z_{r_{i}}\rightarrow Z_{r_{i}-1}\rightarrow ...\rightarrow Z_{r_{i-1}}$ (resp. $Z^{'}_{r_{i}}\rightarrow Z^{'}_{r_{i}-1}\rightarrow ...\rightarrow Z^{'}_{r_{i-1}}$), with $r_{i}>r_{i-1}$ (resp $r_{i}^{'}>r_{i-1}^{'}$), is a brick blow-up for all $i=1, \ldots, t$ as in Definition \ref{DefBriBl}, and the following diagram
\begin{equation*}
\xymatrix{
\widetilde{Z}_{l}\ar[r]\ar[d]_{b} & \widetilde{Z}_{r_{t-1}}\ar[r]\ar[d]_{b_{t-1}} & \widetilde{Z}_{r_{t-2}}\ar[r]\ar[d]_{b_{t-2}} & \cdot\ar[r]\ar[d] & \widetilde{Z}_{r_{1}}\ar[r]\ar[d] & \widetilde{Z}_{0}\ar[d]_{a} \\
\widetilde{Z}^{'}_{l^{'}}\ar[r] & \widetilde{Z}^{'}_{r_{t-1}^{'}}\ar[r] & \widetilde{Z}^{'}_{r_{t-2}^{'}}\ar[r] & \cdot\ar[r] & \widetilde{Z}^{'}_{r_{1}^{'}}\ar[r] & \widetilde{Z}^{'}_{0}
}
\end{equation*}
is commutative.
\end{definition}

\begin{remark}\label{RemAlEqSeqAlEqSM}
If two algebraically marked sequences of blow ups, $(Z_{s},...,Z_{0},\sqcup_{i=1}^{l}F_{i},\pi)_{alg}$ and $(Z^{'}_{s},...,Z^{'}_{0},\sqcup_{i=1}^{l}F^{'}_{i},\pi^{'})_{alg}$, are algebraically equivalent, then their associated algebraically marked sequential morphisms are also algebraically equivalent as in Definition \ref{DefAlEqSM}.
\end{remark}


Roughly speaking, given a partition of the exceptional divisor, combinatorial compatibility with a sequential morphism (Definition \ref{CombCompParSM}) will mean compatibility of the $d-$ary multilinear intersection form, whereas combinatorial compatibility with a sequence of point blow-ups (Definition \ref{CombCompParSeq}) will mean compatibility of proximity relations and degrees of the residue field extensions.\
Before defining formally a combinatorial compatible partition, we show the underlying idea of this concept. Consider a $d-$dimensional smooth algebraic variety
$X$, where $k$ is a perfect field not algebraically closed. Let $P=P_{1}$ a $k-$rational point on $X$ and let $\pi_{1}: Z_{1}\rightarrow X$ be the
blow-up with center $P$. Take a point $P_{2}$ on $Z_{1}$ such that the degree of the field extension $k\subset k(P_{2})\subset k$ is $2$. Let $Z_{1}^{'}$ be the $d-$dimensional smooth algebraic variety defined by $Z_{1}^{'}=Z_{1}\times_{Spec(k)} Spec(K(P_{2}))$. Consider now the blow-up of $Z_{1}$ at $P_{2}$, that is, $\pi_{2}: Z_{2}\rightarrow Z_{1}$. Then, as a consequence of base change, there exists a morphism $\pi_{2}^{'}: Z_{2}^{'}\rightarrow Z_{1}^{'}$, such that $\pi_{2}^{'}$ is the blow-up of $Z_{1}^{'}$ at $P_{2}\times_{Spec(k)} Z_{1}^{'}=\left\{Q_{j}\right\}_{j=2}^{\left[k(P_{2}):k\right]+1}$ and the following diagram is commutative:
\begin{equation*}
\xymatrix{Z_{2}^{'}\ar[r]^{\psi_{2}}\ar[d]_{\pi_{2}^{'}} & Z_{2}\ar[d]^{\pi_{2}} \\
          Z_{1}^{'}\ar[r]^{\psi_{1}} & Z_{1}}
\end{equation*}
Let $H_{1}^{2}, H_{2}^{2}\in Z_{2}$ be the strict transforms of the irreducible exceptional components associated to the sequence $Z_{2}\xrightarrow{\pi_{2}} Z_{1}\xrightarrow{\pi_{1}} Z_{0}=X$, and $H_{1}^{2'}, \left\{H_{j}^{2'}\right\}_{j=2}^{\left[k(P_{2}):k\right]+1}$ the strict transforms of the irreducible exceptional components associated to the sequence $Z_{2}^{'}\xrightarrow{\pi_{2}^{'}} Z_{1}^{'}\xrightarrow{\pi_{1}^{'}} Z_{0}=X$. Then we have that the following properties are satisfied: 
\begin{enumerate}
\item $H_{2}\rightarrow H_{1}$ and  $H_{j}^{'}\rightarrow H_{1}^{'}$ for all $j=2, \ldots, \left[k(P_{2}):k\right]+1$;
\item $(H_{2}^{*})^{d}=\left[k(P_{2}):k\right]=\sum_{j=2}^{\left[k(P_{2}):k\right]+1} (H_{j}^{'*})^{d}$;
\item $\mathcal{I}_{Z_{2}^{'}, E_{Z_{2}^{'}, \pi^{'}}, K}(\overbrace{H_{1}^{'},..., H_{1}^{'}}^{c}, \overbrace{H_{j_{1}}^{'}, \ldots, H_{j_{1}}^{'}}^{d-c})=\mathcal{I}_{Z_{2}^{'}, E_{Z_{2}^{'}, \pi^{'}}, K}(\overbrace{H_{1}^{'},..., H_{1}^{'}}^{c}, \overbrace{H_{j_{2}}^{'}, \ldots, H_{j_{2}}^{'}}^{d-c})$ for all $j_{1}, j_{2}\in\left\{2, \ldots, \left[k(P_{2}):k\right]+1\right\}$;
\item and $H_{j_{1}}^{'}\cap H_{j_{2}}^{'}=\emptyset$ for all $j_{1}, j_{2}\in\left\{2, \ldots, \left[k(P_{2}):k\right]+1\right\}$.
\end{enumerate}
Then, given a sequence of point blow-ups and its associated sequential morphism, it makes sense to introduce the following partitions of the exceptional divisor of combinatorial nature.
\begin{definition}\label{CombCompParSM}
Given a sequential morphism $\pi: Z\rightarrow Z_{0}$ as in Definition \ref{DefSM}, its associated $d-$ary intersection form $\mathcal{I}_{Z, E_{Z, \pi}, K}$, and a partition of the exceptional divisosr $E$, $\sqcup_{i=1}^{l}F_{i}$, we will say that the partition is combinatorially compatible with $\pi$ if for each $i=1, \ldots,l$, and $H_{j_{1}}, H_{j_{2}}\in F_{i}$ there exists a $\sigma\in S_{m}$ such that:
\begin{enumerate}
\item $\sigma(j_{1})=j_{2}$,
\item $\mathcal{I}_{Z, E_{Z, \pi}, K}(H_{i_{1}}, H_{i_{2}},..., H_{i_{d}})=\mathcal{I}_{Z, E_{Z, \pi}, K}(H_{\sigma(i_{1})}, H_{\sigma(i_{2})},..., H_{\sigma(i_{d})})$ for all $i_{1},..,i_{d}$
\end{enumerate}
\end{definition}
Notice that if a given partition is algebraically compatible with a sequential morphism $\pi$, then the partition is also combinatorially compatible with $\pi$.

A combinatorially marked sequential morphism is defined by a pair $(\pi: Z_{s}\rightarrow Z_{0}, \sqcup_{i=1}^{l}F_{i})_{comb}$, where $ \sqcup_{i=1}^{l}F_{i}$ is a partition combinatorially compatible with the sequential morphism $\pi$.

\begin{definition}\label{MInFPar}
Given a combinatorially marked sequential morphism $(\pi: Z_{s}\rightarrow Z_{0}, \sqcup_{i=1}^{l}F_{i})_{comb}$, we can also consider the $d$-ary multilinear intersection form associated to the partition:
\begin{equation*}
\mathcal{I}_{Z, \left(\sqcup_{i=1}^{l}F_{i}\right)_{Z, \pi}, K}: \overbrace{F_{Z}\times F_{Z}\times\cdots\times F_{Z}}^{d}\rightarrow \mathbb{Z},
\end{equation*}
where $F_{Z}$ is the free abelian group generated by $\left\{F_{i}\right\}$, and by an abuse of notation $F_{i}=\sum_{H\in F_{i}} H$. The intersection form is defined by intersecting cycles in the sky $Z$ and taking degrees,
\begin{equation*}
\mathcal{I}_{Z, \left(\sqcup_{i=1}^{l}F_{i}\right)_{Z, \pi}, K}(F_{i_{1}}, F_{i_{2}},..., F_{i_{d}})=deg(F_{i_{1}}\cdot F_{i_{2}}\cdot F_{i_{3}}\cdots F_{i_{d}}),
\end{equation*}
where $F_{i_{1}}\cdot F_{i_{2}}\cdot F_{i_{3}}\cdots F_{i_{d}}$ is a intersection class of $0-$cycles in the abelian group $A_{0}(Z)$, and $deg$ stands for the degree.
\end{definition}

\begin{remark}\label{RemMMap}
Notice that it also makes sense to define the multilinear map associated to a combinatorially marked sequential morphism, and we denote it by $\Phi_{Z, \left(\sqcup_{i=1}^{l}F_{i}\right)_{Z, \pi}, K}$. 
\end{remark}

This fact induces the following equivalence relation.
\begin{definition}\label{ComEqSM}
Given two combinatorially marked sequential morphisms $(\pi: Z_{s}\rightarrow Z_{0}, \sqcup_{i=1}^{l}F_{i})_{comb}$ and $(\pi^{'}: Z_{s}^{'}\rightarrow Z_{0}^{'}, \sqcup_{i=1}^{l}F^{'}_{i})_{comb}$ we say that the associated multilinear maps $\Phi_{Z, \left(\sqcup_{i=1}^{l}F_{i}\right)_{Z, \pi}, K}$ and $\Phi_{Z^{'}, \left(\sqcup_{i=1}^{l}F_{i}^{'}\right)_{Z^{'}, \pi^{'}}, K}$ are equivalent, and we denote it by $\Phi_{Z, \left(\sqcup_{i=1}^{l}F_{i}\right)_{Z, \pi}, K}\sim\Phi_{Z^{'}, \left(\sqcup_{i=1}^{l}F_{i}^{'}\right)_{Z^{'}, \pi^{'}}, K}$, if there exists a $\tau\in\mathcal{S}_{l}$ such that:
\begin{equation*}
\tau(\Phi_{Z, \left(\sqcup_{i=1}^{l}F_{i}\right)_{Z, \pi}, K})=\Phi_{Z^{'}, \left(\sqcup_{i=1}^{l}F_{i}^{'}\right)_{Z^{'}, \pi^{'}}, K}.
\end{equation*}
Moreover, the combinatorially marked sequential morphisms $(\pi: Z_{s}\rightarrow Z_{0}, \sqcup_{i=1}^{l}F_{i})_{comb}$ and $(\pi^{'}: Z_{s}^{'}\rightarrow Z_{0}^{'}, \sqcup_{i=1}^{l}F^{'}_{i})_{comb}$ are said to be combinatorially equivalent, when their associated multilinear maps $\Phi_{Z, \left(\sqcup_{i=1}^{l}F_{i}\right)_{Z, \pi}, K}$ and $\Phi_{Z^{'}, \left(\sqcup_{i=1}^{l}F_{i}^{'}\right)_{Z^{'}, \pi^{'}}, K}$ are equivalent.
\end{definition}

\begin{remark}\label{AlEqSMCombEqMInF}
If two algebraically marked sequential morphisms  $(\pi: Z_{s}\rightarrow Z_{0}, \sqcup_{i=1}^{l}F_{i})_{alg}$ and $(\pi^{'}: Z_{s}^{'}\rightarrow Z_{0}^{'}, \sqcup_{i=1}^{l}F^{'}_{i})_{alg}$ are algebraically equivalent, then by Definition \ref{DefAlEqSM} there exists an isomorphism $b: Z_{s}\rightarrow Z_{s}^{'}$, such that $b(H_{i})=H_{\sigma(i)}^{'}$ for some permutation $\sigma\in S_{m}$. As a result, the class of the associated multilinear map is an invariant. However, the converse is not true. For instance, for $d=2$, let us consider sequences of five point blow-ups, the first one on a rational point of a smooth surface and the other at four different rational points of the exceptional divisor created by the blow-up of the original point. Then, the associated multilinear form, up to a permutation of $S_{5}$, is independent on the choice of the four points; however, two choices with a different cross-ratio provide sequential morphism which are not algebraically isomorphic.
\end{remark}

\begin{definition}\label{CombCompParSeq}
Given a sequence of point blow-ups $(Z_{0},...,Z_{s},\pi)$, its associated weighted directed graph $T_{(Z_{s},...,Z_{0}, \pi)}$, and a partition $\sqcup_{i=1}^{l}F_{i}$, we say that the partition is combinatorially compatible with the sequence $(Z_{0},...,Z_{s},\pi)$ if for each $i=1 \ldots,l$ and $H_{j_{1}}, H_{j_{2}}\in F_{i}$ there exists a $\sigma\in S_{m}$ such that
\begin{enumerate}
\item $\sigma(j_{1})=j_{2}$,
\item $deg(H_{j_{1}})=\left[K(P_{j_{1}}):K\right]=\left[K(P_{\sigma(j_{2})}):K\right]=deg(H_{\sigma(j_{2})})$,
\item if $H_{j_{1}}\in F_{i_{1}}$, $H_{j_{k}}\in F_{i_{k}}$ and $H_{j_{k}}\rightarrow H_{j_{1}}$, then $H_{\sigma(j_{k})}\rightarrow H_{\sigma(j_{1})}$.
\end{enumerate}
\end{definition}
Notice that if a partition is algebraically compatible with a sequence of point blow-ups, then the partition is also combinatorially compatible with the sequence.

A combinatorially marked sequence of point blow-ups is defined by a pair $(Z_{s},...,Z_{0},\sqcup_{i=1}^{l}F_{i},\pi)_{comb}$, where $ \sqcup_{i=1}^{l}F_{i}$ is a partition combinatorially compatible with $(Z_{s},...,Z_{0},\pi)$.

\begin{remark}
Notice that it makes sense to define $F_{i}\rightarrow F_{j}$ if $\exists H_{i}\in F_{i}$, $H_{j}\in F_{j}$ with $H_{i}\rightarrow H_{j}$.
\end{remark}

\begin{definition}
Given a combinatorially marked sequence of point blow-ups, $(Z_{s},...,Z_{0},\sqcup_{i=1}^{l}F_{i},\pi)_{comb}$, it makes sense to consider the weighted directed graph associated to its combinatorially compatible partition. The vertices of $T_{(Z_{s},...,Z_{0}, \sqcup_{i=1}^{l}F_{i}, \pi)}$ are in one to one correspondence with the elements of the partition $F_{1}, F_{2}, \ldots, F_{l}$, and the edges with the pairs $(F_{i}, F_{j})$ such that $j\geq i+1$ and $F_{j}$ is proximate to $F_{i}$. For each $F_{i}$, we denote by $\mathit{l}_{F_{i}}$ the sum of the degrees of the finite field extensions $K\subset K(P_{i_{j}})$, where $K(P_{i_{j}})$ is the residue field of $P_{i_{j}}$, for all $H_{i_{j}}\in F_{i}$.  Then, we define the weighted directed graph associated to the combinatorially marked sequence $(Z_{s},...,Z_{0},\sqcup_{i=1}^{l}F_{i},\pi)_{comb}$  as the directed graph with proximity relations together with the following weights on its vertices: the vertex corresponding to the point $F_{i}$ has the value $\mathit{l}_{F_{i}}$ as weight.
\end{definition}

\begin{definition}\label{Def3}
Given two combinatorially marked sequences of point blow ups, $(Z_{s},...,Z_{0},\sqcup_{i=1}^{l}F_{i},\pi)_{comb}$ and $(Z_{s}^{'},...,Z_{0}^{'},\sqcup_{i=1}^{l}F^{'}_{i},\pi^{'})_{comb}$, with $l=l^{'}$, we say that the associated weighted directed graphs $T_{(Z_{s},...,Z_{0}, \sqcup_{i=1}^{l}F_{i}, \pi)}$, $T_{(Z_{s}^{'},...,Z_{0}^{'}, \sqcup_{i=1}^{l^{'}}F_{i}^{'}, \pi^{'})}$ are equivalent, and we denote it by $T_{(Z_{s},...,Z_{0}, \sqcup_{i=1}^{l}F_{i}, \pi)}\sim T_{(Z_{s}^{'},...,Z_{0}^{'}, \sqcup_{i=1}^{l^{'}}F_{i}^{'}, \pi^{'})}$ , if there exists a $\tau\in\mathcal{S}_{l}$ such that for every two different indexes $i,j$ one has
\begin{enumerate} 
\item there exists a directed edge $F_{i}F_{j}$ if there exists a directed edge $F^{'}_{\tau(i)}F^{'}_{\tau(j)}$,
\item the weight of the vertex associated to $F_{i}$ and that of the vertex associated to $F_{\tau(i)}^{'}$ are the same, that is, $deg(F_{i})=\sum_{H\in F_{i}}deg(H)=\sum_{H^{'}\in F_{i}^{'}}deg(H^{'})=deg(F_{\tau(i)}^{'})$.
\end{enumerate}
 Moreover, we say that two combinatorially marked sequences of point blow-ups $(Z_{s},...,Z_{0},\sqcup_{i=1}^{l}F_{i},\pi)_{comb}$ and $(Z_{s}^{'},...,Z_{0}^{'},\sqcup_{i=1}^{l}F^{'}_{i},\pi^{'})_{comb}$ are combinatorially equivalent when their associated weighted directed graphs $T_{(Z_{s},...,Z_{0}, \sqcup_{i=1}^{l}F_{i}, \pi)}$ and $T_{(Z_{s}^{'},...,Z_{0}^{'}, \sqcup_{i=1}^{l^{'}}F_{i}^{'}, \pi^{'})}$ are equivalent.
\end{definition}

\begin{remark}\label{AlEqSeqCombEqWDG}
If two algebraically marked sequences of point blow-ups $(Z_{s},...,Z_{0},\sqcup_{i=1}^{l}F_{i},\pi)_{alg}$ and $(Z_{s}^{'},...,Z_{0}^{'},\sqcup_{i=1}^{l}F^{'}_{i},\pi^{'})_{alg}$ are algebraically equivalent, then the class of the associated weighted directed graph is an invariant. However, analogously to the case of algebraically marked sequential morphisms and their associated multilinear maps (see Remark \ref{AlEqSMCombEqMInF}), the converse is not true.
\end{remark}

\section{Main results}

Throughout this paper we are restricting ourselves to the case of sequences of point blow-ups over a perfect field $K$, that is sequences of blow-up in which all centers are closed points. Notice that in this case, for each exceptional component $H_{i}$, the value of $r(i)$ (see Definition \ref{DefProx}) is the greatest integer such that the image of $E_{i}$ at $Z_{r(i)}$ is $0-$dimensional.






In order to prove our results, firstly we need to define and characterize the notion of final component of the exceptional divisor $E$. The naive idea is that, given a sequence of blow-ups $(Z_{s},...,Z_{0},\pi)$, an irreducible component $H_{i}$ will be final if there is some other sequence of point blow-ups $(Z_{s}^{'},...,Z_{0}^{'},\pi^{'})$ associated to the sequential morphism $(\pi: Z_{s}\rightarrow Z_{0})$ such that $H_{i}$ will be the exceptional divisor of the last blow-up of $(Z_{s}^{'},...,Z_{0}^{'},\pi^{'})$.

\setcounter{theorem}{0}

\begin{definition}\label{Def7}
Let $(Z_{s},\ldots,Z_{0},\pi)$ be a sequence of blow-ups over $K$ as in Definition \ref{Def1}. The components of the exceptional divisor $E$ in $Z_{s}$ are $\left\{E_{1},...,E_{s}\right\}$. Assume that $H_{i}$ is an irreducible component. Set $H_{i}^{i}$ to be the image of $H_{i}$ in $Z_{i}$. We say that $H_{i}$ is final with respect to $(Z_{s},\ldots,Z_{0},\pi)$ if there exists an open set $U_{i}$ on $Z_{i}$ such that $H_{i}^{i}\subset U_{i}$, $V_{i}=\pi_{s,i}^{-1}(U_{i})\subset Z_{s}$, and $\pi_{s,i}\vert_{V_{i}}: V_{i}\rightarrow U_{i}$ is an isomorphism (see Remark \ref{NoPi} for $\pi_{s,i}$).
\end{definition}

\begin{definition}\label{Def8}
Let $\pi: Z_{s}\rightarrow Z_{0}$ be a sequential morphism. We say that an irreducible component $H_{i}$ of $E$ is final if there exists a sequence of blow-ups $(Z_{s},...,Z_{0},\pi)$ associated to $\pi: Z_{s}\rightarrow Z_{0}$ such that $H_{i}$ is final with respect to this sequence. 
\end{definition}

\begin{lemma}\label{Lem3}
In the particular case of sequences of point blow-ups, if $H_{i}\in E_{i}^{s}$ and $H_{j}\in E_{j}^{s}$ are both final then $H_{i}\cap H_{j}=\emptyset$, that is $H_{i}$ and $H_{j}$ have not geometric points in common over $K$.
\end{lemma}

\begin{proof}
Set $P{i}\in Z_{i-1}$, $P_{j}\in Z_{j-1}$, to be the points such that $H_{i}$ maps to $P_{i}$ and $H_{j}$ maps to $P_{j}$. If $H_{\beta}$ is final with $\beta\in\left\{i,j\right\}$, then $H_{\beta}\cong\mathbb{P}^{d-1}_{K(P_{\beta})}$ and $N_{H_{\beta_{j}}/Z}\cong\mathcal{O}_{H_{\beta}}(-1)$. \\
Let us suppose that $H_{i}\cap H_{j}\neq\emptyset$. Then either $P_{i}$ is proximate to $P_{j}$ or $P_{j}$ is proximate to $P_{i}$. In the first case $H_{j}\not\cong\mathbb{P}^{d-1}_{K(P_{j})}$ and $N_{H_{j_{j}}/Z}\not\cong\mathcal{O}_{H_{j_{j}}}(-1)$ so $H_{j}$ is not final, whereas in the second case $H_{i}\not\cong\mathbb{P}^{d-1}_{K(P_{i})}$ and $N_{H_{i_{j}}/Z}\not\cong\mathcal{O}_{H_{i_{j}}}(-1)$ so $H_{i}$ is not final. Through any of them we get to a contradiction.
\end{proof}

The above result makes a huge difference with respect to the more general case (when considering higher dimensional centers) where two final divisors may have non-empty intersection (see \cite{Camazon25}).

\begin{remark}\label{Rem2}
Assume that we have a sequential morphism $\pi$ associated to a sequence of point blow-ups $(Z_{s},...,Z_{0},\pi)$. If an irreducible component $H_{\alpha}$ of $E$ is final with respect to one representative of the sequences associated to this sequential morphism, then it is final with respect to all. This fact follows directly from Lemma \ref{Lem3}.  
\end{remark}

Before characterizing numerically final divisors, we need a numerical characterization of empty intersections $H_{i}\cap H_{j}=\emptyset$.
\begin{lemma}\label{Lem2}
Let $(Z_{s},\ldots,Z_{0},\pi)$ be a sequence of point blow-ups over $K$. Then, $H_{i}\cap H_{j}=\emptyset$ if and only if $(H_{i})^{s}(H_{j})^{r}=0$ for all $r\neq 0$ and $s\neq 0$, with $r+s=d$.
\end{lemma}

\begin{proof}
If $H_{i}\cap H_{j}=\emptyset$ then $(H_{i})^{s}(H_{j})^{r}=0$ follows directly. \\
In order to prove the necessary condition, we will prove that $H_{i}\cap H_{j}\neq\emptyset\Rightarrow \exists r\neq 0,s\neq 0$ such that $(H_{i})^{s}(H_{j})^{r}\neq 0$, and it is enough to prove it in the case of a sequence of point blow-ups of length $m=3$ since the general result follows by induction.\\
First let $\pi_{1}: Z_{1}\rightarrow Z_{0}$ be the blow-up with center $P_{1}$. Now we blow-up $Z_{1}$ with center $P_{2}$ such that $P_{2}\in E_{1}^{1}$, that is such that $P_{2}$ is proximate to $P_{1}$. Thus we have the following diagram
\begin{equation*}
\xymatrix{
H_{1}^{2}\cap H_{2}^{2}\ar[rr]^{i_{H_{1}^{2}\cap H_{2}^{2},H_{1}^{2}}}\ar[dd] & & H_{1}^{2}\ar[dd]_{\pi_{2}\vert_{H_{1}^{1}}} \\
 & & \\
P_{2}\ar[rr]^{i_{P_{2},H_{1}^{1}}} & & H_{1}^{1}
}
\end{equation*}
Then it follows by Propositions \ref{Pro1} and \ref{Pro6} that
\begin{align*}
H_{1}^{2}\cdot H_{2}^{2}\cdot H_{2}^{2} & =i_{H_{2}^{2},Z_{2} *}(c_{1}(N_{H_{2}^{2}/Z_{2}})\cap H_{1}^{2}\cap H_{2}^{2}), \\
H_{1}^{2}\cdot H_{2}^{2}\cdot H_{1}^{2} & =i_{H_{1}^{2},Z_{2} *}(c_{1}(N_{H_{1}^{2}/Z_{2}})\cap H_{1}^{2}\cap H_{2}^{2}).
\end{align*}
By Proposition \ref{Pro4}
\begin{equation*}
N_{H_{1}^{2}/Z_{2}}\cong\pi_{2}^{*}\vert_{H_{1}^{1}}(N_{H_{1}^{1}/Z_{1}})\otimes\mathcal{O}(-H_{1}^{2}\cap H_{2}^{2}).
\end{equation*}
Furthermore we have the following commutative diagram where all the morphisms are regular embeddings
\begin{equation*}
\xymatrix{
 & H_{1}^{2}\ar[rd]^{j_{H_{1}^{2}, Z_{2}}} & \\
H_{1}^{2}\cap H_{2}^{2}\ar[rr]^{i_{H_{1}^{2}\cap H_{2}^{2}, Z_{2}}}\ar[ru]^{i_{H_{1}^{2}\cap H_{2}^{2}, H_{1}^{2}}}\ar[rd]_{i_{H_{1}^{2}\cap H_{2}^{2}, H_{2}^{2}}} & & Z_{2} \\
 & H_{2}^{2}\ar[ru]_{j_{H_{2}^{2}, Z_{2}}} &
}
\end{equation*}
Now due to the commutativity of the diagram above and the result of Lemma \ref{Lem1} then
\begin{equation*}
H_{1}^{2}\cdot H_{2}^{2}\cdot H_{2}^{2}=i_{H_{1}^{2},Z_{2} *}(c_{1}(N_{H_{1}^{2}\cap H_{2}^{2}/H_{1}^{2}})\cap H_{1}^{2}\cap H_{2}^{2}),
\end{equation*}
so
\begin{align*}
H_{1}^{2}\cdot H_{2}^{2}\cdot H_{2}^{2} & =i_{H_{1}^{2},Z_{2} *}((j_{H_{1}^{2}\cap H_{2}^{2}, H_{1}^{2}*}(H_{1}^{2}\cap H_{2}^{2}))^{2}), \\
H_{1}^{2}\cdot H_{2}^{2}\cdot H_{1}^{2} & =i_{H_{1}^{2},Z_{2} *}(-(j_{H_{1}^{2}\cap H_{2}^{2}, H_{1}^{2}*}(H_{1}^{2}\cap H_{2}^{2}))^{2}).
\end{align*}
By induction on $r$ and $s$ respectively it follows
\begin{align*}
H_{1}^{2}(H_{2}^{2})^{r} & =i_{H_{1}^{2},Z_{2} *}((j_{H_{1}^{2}\cap H_{2}^{2}, H_{1}^{2}*}(H_{1}^{2}\cap H_{2}^{2}))^{r}), \\
(H_{1}^{2})^{s}H_{2}^{2} & =(-1)^{s-1}i_{H_{1}^{2},Z_{2} *}((j_{H_{1}^{2}\cap H_{2}^{2}, H_{1}^{2}*}(H_{1}^{2}\cap H_{2}^{2}))^{2}).
\end{align*}
Finally we can conclude that
\begin{equation*}
(H_{1}^{2})^{s}(H_{2}^{2})^{r}=(-1)^{s-1}deg(i_{H_{1}^{2},Z_{2} *}((j_{H_{1}^{2}\cap H_{2}^{2}, H_{1}^{2}*}(H_{1}^{2}\cap H_{2}^{2}))^{r+s-1}))
\end{equation*}
If we denote $\Delta_{1,2}= deg(i_{H_{1}^{2},Z_{2} *}((j_{H_{1}^{2}\cap H_{2}^{2}, H_{1}^{2}*}(H_{1}^{2}\cap H_{2}^{2}))^{r+s-1}))$, then $(H_{i}^{j})^{s}(H_{j}^{j})^{r}=(-1)^{s-1}\Delta_{1,2}$ for $r+s=d$. Notice that $(H_{1}^{2})^{s}(H_{2}^{2})^{r}\neq 0$. \\
Let $\pi_{3}: Z_{3}\rightarrow Z_{2}$ be the blow-up of $Z_{2}$ with center $P_{3}$, such that $P_{3}\in H_{1}^{2}\cap H_{2}^{2}$, that is $P_{3}$ is proximate to $P_{1}$ and to $P_{2}$. Then it follows that by Theorem \ref{ThmTotTrans}
\begin{equation*}
(H_{1}^{3})^{s}(H_{2}^{3})^{r}=(\pi_{3}^{*}(H_{1}^{2})-j_{H_{3}^{3},Z_{3}*}(H_{3}^{3}))^{s}(\pi_{3}^{*}(H_{2}^{2})-j_{H_{3}^{3},Z*}(H_{3}^{3}))^{r},
\end{equation*}
and due to the projection formula (Propositions \ref{Pro2} and \ref{Pro3}), then
\begin{multline}\label{Eq3P}
(H_{1}^{3})^{s}(H_{2}^{3})^{r}=(\pi_{3}^{*}(H_{1}^{2}))^{s}(\pi_{3}^{*}(H_{2}^{2}))^{r}+(-1)^{d}(j_{H_{3}^{3},Z_{3}*}(H_{3}^{3}))^{d}, \\ =(H_{1}^{2})^{s}(H_{2}^{2})^{r}+(-1)^{d}(j_{H_{3}^{3},Z_{3}*}(H_{3}^{3}))^{d}.
\end{multline}
Since $(H_{1}^{2})^{s}(H_{2}^{2})^{r}\neq0$ and furthermore that it is of the form $(H_{1}^{2})^{s}(H_{2}^{2})^{r}=(-1)^{s-1}\Delta_{1,2}$, then it must exist $r,s$ with $r\neq 0$ and $s\neq 0$ such that $(H_{1}^{3})^{s}(H_{2}^{3})^{r}\neq 0$.\\
For the more general case, let us suppose that $\left\{P_{\alpha_{1}}, P_{\alpha_{2}},..., P_{\alpha_{k}}\right\}$ are proximate to both $P_{1}$ and $P_{2}$. Then by iterating equation (\ref{Eq3P})
\begin{equation*}
(H_{1}^{\alpha_{k}})^{s}(H_{2}^{\alpha_{k}})^{r}=(H_{1}^{2})^{s}(H_{2}^{2})^{r}+(-1)^{d}\sum_{j=1}^{k}(j_{H_{\alpha_{j}}^{\alpha_{j}},Z_{\alpha_{j}}*}(H_{\alpha_{j}}^{\alpha_{j}}))^{d},
\end{equation*}
so it must exist $r,s$ with $r\neq 0$ and $s\neq 0$ such that $(H_{1}^{\alpha_{k}})^{s}(H_{2}^{\alpha_{k}})^{r}\neq 0$.
\end{proof}

Now we are ready to characterize numerically when an irreducible component $H_{i}$ of the exceptional divisor $E$ is final.
\begin{proposition}\label{Pro7}
Let $(Z_{s},\ldots,Z_{0}, \pi)$ (resp. $\pi: Z_{s}\rightarrow Z_{0}$) be a sequence of point blow-ups (resp. sequential morphism) over a perfect field $K$ and let $H_{i}\in E_{i}^{s}$ be an irreducible component of the exceptional divisor of $\pi$. Then $H_{i}$ is final if and only if
\begin{equation*}
(H_{i})^{d}=(-1)^{r}(H_{i})^{s}\cdot (H_{j})^{r}\enspace \text{and}\enspace (H_{i})\cdot(H_{j})^{d-1}>0,
\end{equation*}
for every $j$ such that $H_{i}\cap H_{j}\neq\emptyset$ (see Lemma \ref{Lem2} for a numerical characterization) and for all natural numbers $r$ and $s$ with $r+s=d$.
\end{proposition}

\begin{proof}
We have the following commutative diagram where we denote by $D_{i,j}$ to the scheme theoretic intersection $H_{i}\cap H_{j}$ and all the morphism are regular embeddings
\begin{equation*}
\xymatrix{
 & H_{i}\ar[rd]^{i_{H_{i}, Z}} & \\
D_{i,j}\ar[rr]^{i_{D_{i,j}, Z}}\ar[ru]^{i_{D_{i,j}, H_{i}}}\ar[rd]_{i_{D_{i,j}, H_{j}}} & & Z \\
 & H_{j}\ar[ru]_{i_{H_{j}, Z}} &
}
\end{equation*}
First let us suppose that $H_{i}$ is final .Then $H_{i}\cong\mathbb{P}^{d-1}$ and by Proposition \ref{Pro5} $N_{H_{i}/Z_{n}}=\mathcal{O}_{H_{i}}(-1)$, so it follows by Proposition \ref{Pro1} that
\begin{equation*}
 (H_{i})^{d}=(-1)^{d-1}deg(i_{H_{i}, Z *}(\varsigma^{d-1})),
\end{equation*}
where $\varsigma=c_{1}(\mathcal{O}_{H_{i}}(1))$. \\
Moreover as we have seen in Proposition \ref{Pro1}
\begin{equation*}
H_{i}\cdot H_{j}=i_{H_{i}, Z *}i_{H_{i}, Z}^{*}(H_{j}),
\end{equation*}
and $D_{i,j}$ determines a projective subbundle on $H_{i}$, so by Proposition \ref{Pro6} $D_{i,j}=\varsigma\in A^{1}(H_{i})$. Since $E$ is simple normal crossing divisor, then by Lemma \ref{Lem1}:
\begin{equation*}
N_{H_{j}/Z}\vert_{D_{i,j}}\cong N_{D_{i,j}/H_{i}}. 
\end{equation*}
It follows by Propositions \ref{Pro1} and \ref{Pro6} that
\begin{align*}
H_{i}\cdot H_{j}\cdot H_{j} & =i_{H_{j},Z *}(c_{1}(N_{H_{j}/Z})\cap D_{i,j}) \\
H_{i}\cdot H_{j}\cdot H_{i} & =i_{H_{i}, Z *}(c_{1}(N_{H_{i}/Z})\cap D_{i,j})=(-1)i_{H_{i}, Z *}(\varsigma^{2})
\end{align*}
Now due to the commutativity of the diagram above and the result of Lemma \ref{Lem1} then
\begin{equation*}
H_{i}\cdot H_{j}\cdot H_{j}=i_{H_{i}, Z *}(c_{1}(N_{D_{i,j}/H_{i}})\cap D_{i,j})=i_{H_{i}, Z *}(\varsigma^{2})
\end{equation*}
By induction on $r$ and $s$ respectively it follows
\begin{align*}
H_{i}(H_{j})^{r} & =i_{H_{i}, Z *}(\varsigma^{r}), \\
(H_{i})^{s}H_{j} & =(-1)^{s-1}i_{H_{i}, Z *}(\varsigma^{s}), 
\end{align*}
so we can conclude that
\begin{equation*}
(-1)^{r}(H_{i})^{s}\cdot (H_{j})^{r}=(-1)^{r+s-1}deg(i_{H_{i},Z*}(\varsigma^{r+s-1})).
\end{equation*}
As a result, it is satisfied that $(H_{i})^{d}=(-1)^{r}(H_{i})^{s}\cdot (H_{j})^{r}$. Moreover, $(H_{i})\cdot (H_{j})^{d-1}=deg(i_{H_{i}, Z *}(\varsigma^{d-1}))>0$.

Now let us suppose that $H_{i}$ is not final. If $P_{\alpha}$ is proximate to $P_{i}$, then we have the following commutative diagram
\begin{equation*}
\xymatrix{
 H_{i}^{\alpha}\cap H_{\alpha}^{\alpha}\ar[rr]^{j_{H_{i}^{\alpha}\cap H_{\alpha}^{\alpha},H_{i}^{\alpha}}}\ar[dd] & & H_{i}^{\alpha}\ar[dd]_{\pi_{\alpha}\vert_{H_{i}^{\alpha-1}}}  \\
 & & \\
 P_{\alpha}\ar[rr]^{i_{P_{\alpha},H_{i}^{\alpha-1}}} & & H_{i}^{\alpha-1}  
}
\end{equation*}
Among all the index satisfying $\alpha\rightarrow i$ there must exist an index $j$ such that $j\rightarrow i$ but that there not exists $k$ with $k\rightarrow i$ and $k\rightarrow j$. Let $j$ be such index. Since $H_{i}$ is not final then by Proposition \ref{Pro4} its normal bundle satisfies
\begin{equation*}
N_{H_{i}/Z}=\pi_{n,i}^{*}\vert_{H_{i}^{i}}(N_{H_{i}^{i}/Z_{i}})\otimes\bigotimes_{\alpha\rightarrow i}\pi_{n,\alpha}^{*}\vert_{H_{i}^{\alpha}}(\mathcal{O}(-H_{i}^{\alpha}\cap H_{\alpha}^{\alpha}))
\end{equation*}
Now, by the projection formula (see Proposition \ref{Pro2}) we have that
\begin{equation*}
deg(i_{H_{i}, Z*}(\pi_{n,i}^{*}\vert_{H_{i}^{i}}(c_{1}(N_{H_{i}^{i}/Z_{i}})^{n_{i}}))\prod_{\alpha\rightarrow i}(\pi_{n,\alpha}^{*}\vert_{H_{i}^{\alpha}}((-1)^{d-1}(j_{H_{i}\cap H_{\alpha}^{\alpha},H_{i}^{\alpha}*}(H_{i}^{j}\cap H_{j}^{j})^{d-1}))))=0
\end{equation*}
with $n_{i}+\sum_{\alpha\rightarrow i}n_{\alpha}=d$, so
\begin{multline*}
(H_{i})^{d}=deg(i_{H_{i},Z*}(\pi_{n,i}^{*}\vert_{H_{i}^{i}}(c_{1}(N_{H_{i}^{i}/Z_{i}})^{d-1}))+ \\ \sum_{\alpha\rightarrow i}deg(i_{H_{i},Z*}(\pi_{n,\alpha}^{*}\vert_{H_{i}^{\alpha}}((-1)^{d-1}(j_{H_{i}\cap H_{\alpha}^{\alpha},H_{i}^{\alpha}*}(H_{i}^{\alpha}\cap H_{\alpha}^{\alpha})^{d-1}))))
\end{multline*}
Furthermore, by an analogous reasoning to that of the case where $H_{i}$ is final, then
\begin{align*}
H_{i}\cdot H_{j}\cdot H_{j} & =i_{H_{i}, Z *}(c_{1}(N_{D_{i,j}/H_{i}})\cap D_{i,j}), \\ & =i_{H_{i}, Z *}((j_{H_{i}^{j}\cap H_{j}^{j}, H_{i}^{j}*}(H_{i}^{j}\cap H_{j}^{j}))^{2}), \\ & =i_{H_{i}, Z *}((j_{H_{i}^{j}\cap H_{j}^{j}, H_{i}^{j}*}(D_{i,j}^{2})); \\
H_{i}\cdot H_{j}\cdot H_{i} & =i_{H_{i}, Z *}(c_{1}(N_{H_{i}/Z})\cap D_{i,j}), \\ & =i_{H_{i}, Z *}(-(j_{H_{i}^{j}\cap H_{j}^{j}, H_{i}^{j}*}(H_{i}^{j}\cap H_{j}^{j}))^{2}), \\ & =i_{H_{i}, Z *}(-(j_{H_{i}^{j}\cap H_{j}^{j}, H_{i}^{j}*}(D_{i,j}))^{2}).
\end{align*}
By induction on $r$ and $s$ respectively it follows
\begin{align*}
H_{i}(H_{j})^{r} & =i_{H_{i}, Z *}((j_{H_{i}^{j}\cap H_{j}^{j}, H_{i}^{j}*}(D_{i,j}^{r})),\\
(H_{i})^{s}H_{j} & =(-1)^{s-1}i_{H_{i}, Z *}((j_{H_{i}^{j}\cap H_{j}^{j}, H_{i}^{j}*}(D_{i,j}))^{s}),
\end{align*}
so we have that
\begin{equation*}
(H_{i})^{s}(H_{j})^{r}=(-1)^{s-1}deg(i_{H_{i}, Z *}(j_{H_{i}^{j}\cap H_{j}^{j}, H_{i}^{j}*}(D_{i,j})^{r+s-1})).
\end{equation*}
If $d$ is even then
\begin{equation*}
(-1)^{r}(H_{i})^{s}(H_{j})^{r}=(-1)^{r+s-1}deg(i_{H_{i}, Z *}(j_{H_{i}^{j}\cap H_{j}^{j}, H_{i}^{j}*}(D_{i,j})^{r+s-1}))\neq(H_{i})^{d}
\end{equation*}
since
\begin{multline*} 
deg(i_{H_{i},Z*}(\pi_{n,i}^{*}\vert_{H_{i}^{i}}(c_{1}(N_{H_{i}^{i}/Z_{i}})^{d-1})))+ \\ \sum_{\substack{\alpha\rightarrow i \\ \alpha\neq j}}deg(i_{H_{i},Z*}(\pi_{n,\alpha}^{*}\vert_{H_{i}^{\alpha}}((-1)^{d-1}(j_{H_{i}\cap H_{\alpha}^{\alpha},H_{i}^{\alpha}*}(H_{i}^{\alpha}\cap H_{\alpha}^{\alpha})^{d-1}))))<0
\end{multline*}
If $d$ is odd then 
\begin{equation*}
(H_{i})\cdot (H_{j})^{d-1}=deg(i_{H_{i}, Z *}((j_{H_{i}^{j}\cap H_{j}^{j}, H_{i}^{j}*}(D_{i,j})^{r+s-1}))),
\end{equation*}
so $(H_{i})(H_{j})^{d-1}<0$.
\end{proof}

\begin{remark}
Notice that whereas in \cite{CampilloReguera94} the authors show two different numerical characterizations for final divisors, depending on parity of $d$, our numerical characterization works for both cases.
\end{remark} 

Next, we prove that our numerical characterization of final divisors can be extended in a natural manner to the elements of an algebraically compatible partition of $E$.
\begin{proposition}
Let $(Z_{s},...,Z_{0},\sqcup_{i=1}^{l}F_{i},\pi)_{alg}$ be an algebraically marked sequence, and let $H, H^{'}\in F_{i}$. Then $H$ is final if and only if $H^{'}$ is also final.
\end{proposition}

\begin{proof}
If $H, H^{'}\in F_{i}$ then there exists a sequence $(\widetilde{Z_{s}},...,\widetilde{Z_{0}},\widetilde{E},\widetilde{\pi})$ over $\widetilde{K}$ such that $\beta(H)=\beta(H^{'})$, where $\beta: \widetilde{Z}_{s}\times_{Spec(\widetilde{K})} Spec(K)\rightarrow\widetilde{Z}_{s}$, so it follows that if $H$ satisfies the numerical condition of Proposition \ref{Pro7}, $H^{'}$ will satisfy it too.
\end{proof}

\begin{proposition}
Let $(Z_{s},...,Z_{0},\sqcup_{i=1}^{l}F_{i},\pi)_{alg}$ be an algebraically marked sequence. Then $F_{i}$ is final if and only if
\begin{equation*}
(F_{i})^{d}=(-1)^{r}(F_{i})^{s}\cdot (F_{j})^{r}\enspace \text{and}\enspace (F_{i})\cdot(F_{j})^{d-1}>0,
\end{equation*}
for every $j$ such that $F_{i}\cap F_{j}\neq\emptyset$ and for all natural numbers $r$ and $s$ with $r+s=d$.
\end{proposition}

\begin{proof}
This is a consequence of Proposition \ref{Pro7}, since if $H, H^{'}\in F_{i}$ and $H\neq H^{'}$ then $H\cap H^{'}=\emptyset$.
\end{proof}

\begin{corollary}\label{Cor1}
Let $(\pi: Z_{s}\rightarrow Z_{0}, \sqcup_{i=1}^{l}F_{i})_{alg}$ and $(\pi^{'}: Z_{s}^{'}\rightarrow Z_{0}^{'}, \sqcup_{i=1}^{l^{'}}F_{i}^{'})_{alg}$ algebraically marked sequential morphisms that are algebraically equivalent, and $b$ the isomorphism $b: Z_{s}\rightarrow Z_{s}^{'}$ (see Definition \ref{DefAlEqSM}). Then an element of the algebraically compatible partition $F_{i}$ is final if and only if $b(F_{i})$ is also final.
\end{corollary}

Now, we define a key tool for our proof of Theorems \ref{Thm1} and \ref{Thm21}, that of a regular projective contraction.
\begin{definition}\label{Def9}
Let $X^{'}$ be a $d-$dimensional variety, $L$ be an effective Cartier divisor on $X^{'}$ and $Y$ an $r-$dimensional variety with $r<d-1$. Then we say that $L$ is contractible to $Y$ within $X^{'}$ if there exists a variety $X$ and a proper birational morphism $f: X^{'}\longrightarrow X$ such that:
\begin{enumerate}
\item $f(L)=Y$; 
\item $L$ is a closed subset of $X^{'}$ which consists of the points where $f$ is not an isomorphism.
\end{enumerate}
We call this triple $(X^{'}, f, X)$ a contraction of $L$ to $Y$, or simply a contraction. We shall say that $L$ is regularly and projectively contractible to $Y$ within $X^{'}$, when moreover 
\begin{enumerate}
\setcounter{enumi}{2}
\item X is a non-singular projective variety.
\end{enumerate}
\end{definition}

Before going on, we need to prove the following technical lemma that is crucial for the uniqueness of the regular projective contractions.
\begin{lemma}\label{Lem4}
Let $X$ and $Y$ be two affine normal varieties such that $X=Spec(A)$ and $Y=Spec(B)$. Let $\pi: Z\rightarrow X$ and $\pi^{'}: Z^{'}\rightarrow Y$ be proper morphisms. If $Z$ is isomorphic to $Z^{'}$ then $A\cong B$.
\end{lemma}

\begin{proof}
By \cite[Theorem 3.2.1]{EGAIII.1}, since $\pi$ and $\pi^{'}$ are proper morphisms then $\pi_{*}(\mathcal{O}_{Z})$ and $\pi_{*}^{'}(\mathcal{O}_{Z^{'}})$ are a coherent sheaves on $X$ and $Y$ respectively. Since $X$ and $Y$ are both normal, then $\mathcal{O}_{X}\cong\pi_{*}(\mathcal{O}_{Z})$ and $\mathcal{O}_{Y}\cong\pi_{*}^{'}(\mathcal{O}_{Z^{'}})$, so 
\begin{gather*}
A\cong H^{0}(X, \mathcal{O}_{X})\cong H^{0}(Z, \mathcal{O}_{Z}), \\ 
B\cong H^{0}(Y, \mathcal{O}_{Y})\cong H^{0}(Z^{'}, \mathcal{O}_{Z^{'}}).
\end{gather*}
Since $H^{0}(Z, \mathcal{O}_{Z})\cong H^{0}(Z^{'}, \mathcal{O}_{Z^{'}})$, then it follows that $A\cong B$.
\end{proof}

We are now ready to prove that by contracting a final divisor, the variety obtained corresponds to the sky of a sequence of point blow-ups with the same ground variety as the original one, at the same time that the image of the exceptional divisor preserves its normal crossing structure.
\begin{proposition}\label{Pro8}
Let $(Z_{s},...,Z_{0}, \pi)$ be a sequence of point blow-ups (as in Definition \ref{Def1}) of length $m$ and let $H_{i}\in E_{i}^{s}$ be an irreducible component of the exceptional divisor of $\pi$. If $H_{i}$ is final, then there exists a regular projective contraction $(Z, f_{m}, X_{m-1})$ of $H_{i}$ to a point such that $f_{m}(E_{Z})$ is a simple normal crossing divisor and $X_{m-1}$ is the sky of a sequence of point blow-ups with ground $Z_{0}$.
\end{proposition}

\begin{proof}
Since $H_{i}$ is final there must exists an isomorphism between the two opens sets $U_{i}\subset Z_{i}$ and $V_{i}\subset Z_{s}$ via $\pi_{s,i}$. After shrinking $U_{i}$ if necessary, we may assume that $U_{i}\setminus H_{i}^{i}$ is isomorphic via $\pi_{i}$ to an open set of $Z_{i-1}\setminus\left\{P_{i}\right\}$ where $P_{i}=\pi_{i}(H_{i}^{i})$.\\
Note that $W_{i}=\pi_{i}(U_{i})$ is an open set in $Z_{i-1}$. In fact $\pi_{i}\vert_{U_{i}}$ is the blow-up of $W_{i}$ at $P_{i}$.  
$$\xymatrix{V_{i}\ar[r]^{\pi_{s,i}\vert_{V_{i}}} & U_{i}\ar[l]\ar[d]_{\pi_{i}\vert_{U_{i}}} \\
 & W_{i}} $$\\
Set $\phi=(\pi_{i}\circ\pi_{s,i})\vert_{V_{i}}$ the composition morphism from $V_{i}$ to $W_{i}$
$$\xymatrix{V_{i}\ar[r]^{\pi_{s,i}\vert_{V_{i}}}\ar[rd]^{\phi} & U_{i}\ar[d]^{\pi_{i}\vert_{U_{i}}}\ar[l] \\
& W_{i}} $$\\
where $\phi:=\pi_{i}\circ\pi_{s,i}$. \\
Set $\overline{W}_{i}=Z\setminus H{i}$. We construct $X_{m-1}$ by gluing $W_{i}$ and $\overline{W}_{i}$ along the open isomorphic sets $W_{i}\setminus\left\{P_{i}\right\}\subset W_{i}$ and $V_{i}\setminus H_{i}\subset\overline{W}_{i}$. Note that $W_{i}\setminus\left\{P_{i}\right\}\stackrel{\pi_{i}}{\cong} U_{i}\setminus H_{i}^{i}\stackrel{\pi_{s,i}}{\cong} V_{i}\setminus H_{i}$.\\
Now we define $f_{m}: Z\rightarrow X_{m-1}$, $f_{m}\vert_{\overline{W}_{i}}=Id_{\overline{W}_{i}}$, $f_{m}\vert_{V_{i}}=\phi$, which is well defined by the isomorphisms. \\
Finally, it is clear from the construction that if we denote by $D_{X_{m-1}}$ to the image $f_{m}(E_{Z})$, then $D_{X_{m-1}}$ is a simple normal crossing divisor.

\underline{\textsl{An alternative construction of the contraction if $K=\overline{K}$}} \\
Since $H_{i}$ is final, then $H_{i}\cong\mathbb{P}_{K(P_{i})}^{d-1}$, where $P_{i}=\pi_{s,i}(H_{i})$, and moreover by Proposition \ref{Pro5} its normal bundle $N_{H_{i}/Z}\cong\mathcal{O}_{H_{i}}(-1)$.
Let $\mathcal{F}$ be a very ample line bundle on $Z$. Then $\mathcal{F}\otimes\mathcal{O}_{H_{Z,i}}=\mathcal{L}_{2}\otimes\mathcal{O}_{H_{Z,i}}(u)$. If we consider the line bundle $\mathcal{L}:=\mathcal{F}\otimes\mathcal{O}(H_{Z,i})^{\otimes u}$, then by \cite[Corollary 2.]{Ishii77} there exists a regular projective contraction $(Z,\varphi, X_{m-1}^{'})$ of $H_{Z,i}$ to a closed point, that we will denote by $P_{i}^{'}$, with $\varphi\vert_{H_{Z,i}}=P_{i}^{'}$, such that $\varphi$ is defined by the complete linear system $\left|\mathcal{L}\right|$. To see that $D_{X_{m-1}^{'}}:=\varphi(E_{Z})$ is still a simple normal crossing divisor we prove that the contraction is unique up to isomorphism. By \cite[Theorem 3]{Ishii77} $(Z, \varphi)$ is the blowing up of $X_{m-1}^{'}$ at a point $P_{i}^{'}$. Let $Y_{i}$ be an affine open neighborhood of $P_{i}$ in $X_{m-1}$. If we denote by $Y_{i}^{'}:= \varphi(f_{m}^{-1}(Y_{i}))$, then is $Y_{i}^{'}$ is an affine neighborhood of $P_{i}^{'}$ in $X_{m-1}^{'}$ and there exist two proper morphisms
$$\xymatrix{ & f_{m}^{-1}(Y_{i})\ar[dl]_{f_{m}\vert_{f_{m}^{-1}(Y_{i})}}\ar[dr]^{\varphi\vert_{f_{m}^{-1}(Y_{i})}} & \\
Y_{i} & & Y_{i}^{'}}$$\\
By Lemma \ref{Lem4}  this implies that $Y_{i}\cong Y_{i}^{'}$, so its then clear that $D_{X_{m-1}^{'}}$ is a simple normal crossing divisor.

So we have proved that there exists a regular projective contraction $(Z, f_{m}, X_{m-1})$ of $H_{i}$ to a point $P_{i}\in X_{m-1}$, $f_{m}: Z_{s}\rightarrow X_{m-1}$.
Following the notations of Definition \ref{Def7}, let $W_{i}=\pi_{i}(U_{i})$. Then by Definition \ref{Def9} $f_{m}\vert_{Z\setminus V_{i}}: Z\setminus V_{i}\rightarrow X_{m-1}\setminus f_{m}(V_{i})$ is an isomorphism. Now we define $g: X_{m-1}\rightarrow Z_{i-1}$ as follows: $g\vert_{\overline{W}_{i}}=\pi_{s,i-1}\vert_{\overline{W}_{i}}$ and $g\vert_{W_{i}}=Id_{W_{i}}$. By our construction of $X_{m-1}$ $g$ is well defined, and by the definition $g: X_{m-1}\rightarrow Z_{i-1}$ is a sequence of point blow-ups.\\
Hence the composition $X_{m-1}\rightarrow Z_{i-1}\rightarrow Z_{0}$ is a sequence of blow-ups.

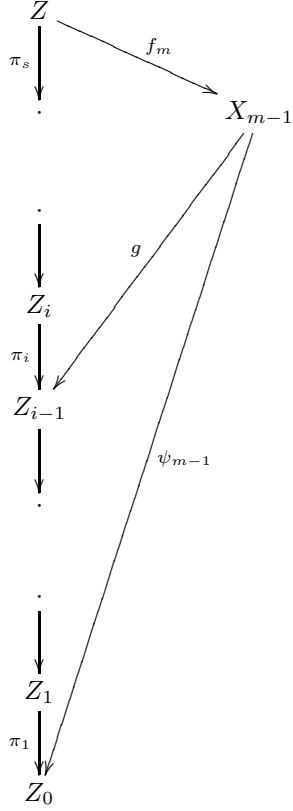
\begin{equation*}
\xymatrix{Z\ar[d]_{\pi_{s}}\ar[rrd]^{f_{m}} & & \\
\cdot & & X_{m-1}\ar[llddddd]^{\psi_{m-1}}\ar[lldd]_{g}\\
\cdot\ar[d] & & \\
Z_{i-1}\ar[d] & & \\
\cdot & & \\
\cdot\ar[d] & & \\
Z_{0} & & }
\end{equation*}
\end{proof}

The previous result joint with Proposition \ref{Pro7} lead us to prove one on the main theorems of the paper: There exists a bijection $\Theta: (\Phi_{Z, \left(\sqcup_{i=1}^{l}F_{i}\right)_{Z, \pi}, K})/\sim\rightarrow (T_{(Z_{s},...,Z_{0}, \sqcup_{i=1}^{l}F_{i}, \pi)})/\sim$ between the equivalence classes of the combinatorial objects associated to combinatorially marked sequences of point blow-ups and the corresponding sequential morphisms.
\begin{theorem}\label{Thm1}
Let $(Z_{s},...,Z_{0},\sqcup_{i=1}^{l}F_{i},\pi)_{comb}$ and $(Z^{'}_{s},...,Z^{'}_{0},\sqcup_{i=1}^{l^{'}}F_{i}^{'},\pi^{'})_{comb}$ be two combinatorially marked sequences of point blow-ups, with $l=l^{'}$. The associated weighted directed graphs $T_{(Z_{s},...,Z_{0}, \sqcup_{i=1}^{l}F_{i}, \pi)}$, $T_{(Z_{s}^{'},...,Z_{0}^{'}, \sqcup_{i=1}^{l^{'}}F_{i}^{'}, \pi^{'})}$ are equivalent if and only if the multilinear maps $\Phi_{Z, \left(\sqcup_{i=1}^{l}F_{i}\right)_{Z, \pi}, K}$ and $\Phi_{Z^{'}, \left(\sqcup_{i=1}^{l}F_{i}^{'}\right)_{Z^{'}, \pi^{'}}, K}$ associated to the corresponding combinatorially marked sequential morphisms are also equivalent as in Definition \ref{ComEqSM}.
\end{theorem}

First we will prove that if the multilinear maps $\Phi_{Z, \left(\sqcup_{i=1}^{l}F_{i}\right)_{Z, \pi}, K}$ and $\Phi_{Z^{'}, \left(\sqcup_{i=1}^{l^{'}}F_{i}^{'}\right)_{Z^{'}, \pi^{'}}, K}$ associated to the corresponding combinatorially marked sequential morphisms $(\pi: Z_{s}\rightarrow Z_{0}, \sqcup_{i=1}^{l}F_{i})_{comb}$ and $(\pi^{'}: Z_{s}^{'}\rightarrow Z_{0}^{'}, \sqcup_{i=1}^{l^{'}}F^{'}_{i})_{comb}$ are equivalent, then the weighted directed graphs $T_{(Z_{s},...,Z_{0}, \sqcup_{i=1}^{l}F_{i}, \pi)}$, $T_{(Z_{s}^{'},...,Z_{0}^{'}, \sqcup_{i=1}^{l^{'}}F_{i}^{'}, \pi^{'})}$ associated to the corresponding combinatorially marked sequences of points blow-ups are equivalent too. To begin with, we need to establish a numerical characterization of proximity.
\begin{lemma}\label{LemNumPro}
Let $(Z_{s},\ldots,Z_{0},\sqcup_{i=1}^{l}F_{i},\pi)_{comb}$ be a combinatorially marked sequence of point blow-ups. Then $P_{i}\rightarrow P_{j}$ if and only if
\begin{enumerate}
\item there exists $\alpha\in\left\{2,3,...,m-1,m\right\}$ such that $H_{X_{\alpha},i}\cap H_{X_{\alpha},j}\neq\emptyset$ (see numerical characterization of Lemma \ref{Lem2}), \label{Cond1}
\item and $(H_{X_{\alpha},i})^{d}=(-1)^{r}(H_{X_{\alpha},i})^{s}\cdot (H_{X_{\alpha},k})^{r}$, $(H_{X_{\alpha},i})\cdot(H_{X_{\alpha},k})^{d-1}>0$, for all $k$ such that $H_{X_{\alpha},i}\cap H_{X_{\alpha},k}\neq\emptyset$, \label{Cond2}
\end{enumerate}
where $Z_{s}=X_{m}\rightarrow X_{m-1}\rightarrow\cdots\rightarrow X_{\alpha}\rightarrow\cdots\rightarrow X_{0}=Z_{0}$ is any of the sequences of regular projective contractions obtained by iterating the process above in Proposition \ref{Pro8}. The sequence is obtained by contracting a final divisor at each step.
\end{lemma}

\begin{proof}
If $P_{i}$ is proximate to $P_{j}$ then $H_{i}^{i}\cap H_{j}^{i}\neq\emptyset$ so \ref{Cond1} holds. Moreover, if $P_{i}\in Z_{r}$, then $H_{i}^{i}$ is final for the sequence $\pi_{r+1}\circ\pi_{r}\circ\cdot\cdot\cdot\circ\pi_{1}$, for some $r\geq i$ and we have \ref{Cond2}. \\
Conversely, if $H_{X_{\alpha},i}$ is final for the sequence $\pi_{r+1}\circ\pi_{r}\circ\cdot\cdot\cdot\circ\pi_{1}$ for some $r\geq i$ then by Proposition \ref{Pro8} there exists a regular projective contraction $f_{\alpha}: X_{\alpha}\rightarrow X_{\alpha-1}$ of $H_{X_{\alpha},i}$ such that $f_{\alpha}(H_{X_{\alpha},i})=O_{i,\alpha-1}\subset H_{X_{\alpha-1},j}$.
\end{proof}

\begin{remark}
Notice that the previous characterization can be naturally adapted in order to characterize numerically the proximity between two elements of a combinatorially compatible partition. Hence, it makes sense to denote $F_{i}\rightarrow F_{j}$.
\end{remark}

\begin{proof}[Proof of Thm. \ref{Thm1}]
Assume that the multilinear maps $\Phi_{Z, \left(\sqcup_{i=1}^{l}F_{i}\right)_{Z, \pi}, K}$ and $\Phi_{Z^{'}, \left(\sqcup_{i=1}^{l}F_{i}^{'}\right)_{Z^{'}, \pi^{'}}, K}$, associated to the corresponding combinatorially marked sequential morphisms $(\pi: Z_{s}\rightarrow Z_{0}, \sqcup_{i=1}^{l}F_{i})_{comb}$ and $(\pi^{'}: Z_{s}^{'}\rightarrow Z_{0}^{'}, \sqcup_{i=1}^{l^{'}}F^{'}_{i})_{comb}$, are combinatorially equivalent. If $F_{i}$ is final, then there exists a $\tau\in S_{l}$ such that:
\begin{enumerate}
\item $F^{'}_{\tau(i)}$ is final,
\item $F_{i}\cap F_{\beta}\neq\emptyset$ if and only if $F^{'}_{\tau(i)}\cap F^{'}_{\tau(\beta)}\neq\emptyset$,
\item  and $F_{\beta_{1}}\cdot F_{\beta_{2}}\cdots F_{\beta_{d}}=F_{\tau(\beta_{1})}\cdot F_{\tau(\beta_{2})}\cdot\cdot\cdot F_{\tau(\beta_{d})}$.
\end{enumerate}
Let $f_{l}: Z\rightarrow X_{l-1}$ ($f_{l}^{'}: Z^{'}\rightarrow X_{l-1}^{'}$) be the contraction of all the irreducible components $H_{\delta}\in F_{i}$ ($H_{\delta}^{'}\in F_{\tau(i)}^{'}$ respectively). Recall that we denote $\psi_{l-1}: X_{l-1}\rightarrow Z_{0}$ ($\psi_{l-1}^{'}: X_{l-1}^{'}\rightarrow Z_{0}^{'}$) to the sequential morphism associated to the sequence of point blow-ups with sky $X_{l-1}$ ($X_{l-1}^{'}$ respectively) and ground variety $Z_{0}$ ($Z_{0}^{'}$ respectively) whose existence is proved in Proposition \ref{Pro8}. As a result of Theorem \ref{ThmTotTrans} and Proposition \ref{Pro3}, and by an abuse of notation we have that:
\begin{align*}
F_{X_{l-1},\beta_{1}}\cdot F_{X_{l-1},\beta_{2}}\cdots F_{X_{l-1},\beta_{d}} & =(F_{\beta_{1}}+\delta_{\beta_{1},i}F_{i})\cdot (F_{\beta_{2}}+\delta_{\beta_{2},i}F_{i})\cdots (F_{\beta_{d}}+\delta_{\beta_{d},i}F_{i}), \\
F_{X_{l-1}^{'},\beta_{1}}\cdot F_{X_{l-1}^{'},\beta_{2}}\cdots F_{X_{l-1}^{'},\beta_{d}} & =(F_{\beta_{1}}^{'}+\delta_{\beta_{1},\tau(i)}^{'}F_{\tau(i)})\cdot (F_{\beta_{2}}^{'}+\delta_{\beta_{2},\tau(i)}^{'}F_{i}^{'})\cdots (F_{\beta_{d}}^{'}+\delta_{\beta_{d},\tau(i)}^{'}F_{i}^{'}),
\end{align*}
where $\delta_{\beta_{j},i}$ ($\delta_{\beta_{j},\tau(i)}^{'}$) is equal to $1$ if $F_{i}\rightarrow F_{\beta_{j}}$ ($F_{\tau(i)}^{'}\rightarrow F_{\beta_{j}}^{'}$ respectively) and $0$ otherwise. Then, it follows that there exists a $\widetilde{\tau}\in S_{l-1}$ such that
\begin{equation*}
F_{X_{l-1},\beta_{1}}\cdot F_{X_{l-1},\beta_{2}}\cdots F_{X_{l-1},\beta_{d}}=F^{'}_{X_{l-1}^{'},\widetilde{\tau}(\beta_{1})}\cdot F^{'}_{X_{l-1}^{'},\widetilde{\tau}(\beta_{2})}\cdots F^{'}_{X_{l-1}^{'},\widetilde{\tau}(\beta_{d})}.
\end{equation*}
By iterating the above process, we obtain that $\Phi_{X_{\alpha},\left(\sqcup_{i=1}^{\alpha}F_{i}^{\alpha}\right)_{X_{\alpha}, \psi_{\alpha}}, K}\sim\Phi_{X_{\alpha}^{'},\left(\sqcup_{i=1}^{\alpha}F_{i}^{\alpha '}\right)_{X_{\alpha}^{'}, \psi_{\alpha}^{'}}, K}$ for $\alpha=1, \ldots, l-1$.Then, as a consequence of Lemma \ref{LemNumPro}, any two equivalent multilinear maps preserve proximity relations. Moreover, it is verified that $deg(F_{i})=deg(F^{'}_{\tau(i)})$, so equivalent multilinear maps preserve also degrees.\\
Conversely, we assume now that the weighted directed trees, $T_{(Z_{s},...,Z_{0}, \sqcup_{i=1}^{l}F_{i}, \pi)}$ and $T_{(Z_{s}^{'},...,Z_{0}^{'}, \sqcup_{i=1}^{l^{'}}F_{i}^{'}, \pi^{'})}$, associated to the corresponding combinatorially marked sequences of points blow-ups are equivalent. We want to prove that the multilinear maps $\Phi_{Z, \left(\sqcup_{i=1}^{l}F_{i}\right)_{Z, \pi}, K}$ and $\Phi_{Z^{'}, \left(\sqcup_{i=1}^{l}F_{i}^{'}\right)_{Z^{'}, \pi^{'}}, K}$,   associated to the corresponding combinatorially marked sequential morphisms $(\pi: Z_{s}\rightarrow Z_{0}, \sqcup_{i=1}^{l}F_{i})_{comb}$ and $(\pi^{'}: Z_{s}^{'}\rightarrow Z_{0}^{'}, \sqcup_{i=1}^{l^{'}}F^{'}_{i})_{comb}$, are also equivalent. Firstly, by and abuse of notation we denote by $F_{j}^{*}$ ($F_{j}^{'*}$) to $H_{j,1}^{*}+\ldots H_{j,n}^{*}$ ($H_{j,1}^{'*}+\ldots H_{j,n^{'}}^{'*}$ respectively) for all $H_{j,k}\in F_{j}$ ($H_{j,k^{'}}^{'}\in F_{j}^{'}$ respectively). Then, there exists a $\sigma\in S_{l}$, such that by applying iteratively Theorem \ref{ThmTotTrans}, it follows that:
\begin{align*}
& F_{i}= F_{i}^{*}-\sum_{\beta\rightarrow i}F_{\beta}^{*}, \\
& F^{'}_{\sigma(i)}= F^{'*}_{\sigma(i)}-\sum_{\sigma(\beta)\rightarrow \sigma(i)}F^{'*}_{\sigma(\beta)}.
\end{align*}
Moreover, as a consequence of Proposition \ref{Pro2}, it holds that $F_{\beta_{1}}^{*}F_{\beta_{2}}^{*}\cdots F_{\beta_{d}}^{*}\neq 0\enspace$ if and only if $\beta_{1}=\beta_{2}=\ldots=\beta_{d}$. Now, if an element of the partition $F_{i}$ is final then $F_{i}=F_{i}^{*}$, so $deg(F_{i}^{*})=deg(F^{'*}_{\sigma(i)})$ for all $i=1, \ldots, l$. Finally, as a result of Theorem \ref{ThmTotTrans} we have that:
\begin{equation*}
F_{\beta_{1}}\cdot F_{\beta_{2}}\cdots F_{\beta_{d}}=(F_{\beta_{1}}^{*}-\sum_{\delta\rightarrow\beta_{1}}F_{\delta}^{*})\cdot (F_{\beta_{2}}^{*}-\sum_{\delta\rightarrow\beta_{2}}F_{\delta}^{*})\cdots (F_{\beta_{d}}^{*}-\sum_{\delta\rightarrow\beta_{d}}F_{\delta}^{*}),
\end{equation*}
so it follows that $F_{\beta_{1}}\cdot F_{\beta_{2}}\cdot\cdot\cdot F_{\beta_{d}}=F_{\sigma(\beta_{1})}^{'}\cdot F_{\sigma(\beta_{2})}^{'}\cdots F_{\tau(\sigma_{d})}^{'}$.
\end{proof}

Now, we go an step further, and the following results leads us to prove Theorem \ref{Thm21}, that is, there exists a thinner bijection between the algebraic equivalence classes of algebraically marked sequences of point blow-ups and that of the associated sequential morphisms.
\begin{proposition}\label{Thm41}
Let $(\pi: Z_{s}\rightarrow Z_{0},\sqcup_{i=1}^{l}F_{i})_{alg}$ be an algebraically marked sequential morphism. Given the $d-$ary multilinear intersection form associated to the partition, $\mathcal{I}_{Z, \left(\sqcup_{i=1}^{l}F_{i}\right)_{Z, \pi}, K}$ (see Definition \ref{MInFPar}), we can recover all the algebraically marked sequences of point blow-ups that are associated to an algebraically marked sequential morphisms in the same algebraic equivalence class of $(\pi: Z_{s}\rightarrow Z_{0},\sqcup_{i=1}^{l}F_{i})_{alg}$.
\end{proposition}

\begin{proof}[Proof of Proposition \ref{Thm41}, Part $I$]
Since $\sqcup_{i=1}^{l}F_{i}$ is a partition algebraically compatible with $\pi$ then $\exists\widetilde{K}\subset K$ as in Definition \ref{AlComPartSM}. If $H\in F_{i}$ is final then $\widetilde{H}=\beta(H)$ is final for $\widetilde{\pi}: \widetilde{Z}_{s}\rightarrow\widetilde{Z}_{0}$. We will prove this result first by contracting one irreducible component of the exceptional divisor $\widetilde{E}$ at a time.\\
Since the set formed by final divisors is not empty, let us suppose that $\widetilde{H}_{i}$ is final, then by Proposition \ref{Pro7} there exists a regular projective contraction $(\widetilde{Z}_{s}, \widetilde{f_{l}}, \widetilde{X}_{l-1})$ of $\widetilde{H}_{i}$ to a point such that $\widetilde{X}_{l-1}$ is the sky of a sequence of point blow-ups with ground $\widetilde{Z}_{0}$ .

\begin{equation*}
\xymatrix{\widetilde{Z}_{s}\ar[d]_{\widetilde{\pi}_{s}}\ar[rrd]^{\widetilde{f}_{l}} & & \\
\cdot & & \widetilde{X}_{l-1}\ar[llddddddd]^{\widetilde{\psi}_{l-1}}\\
\cdot\ar[d] & & \\
\widetilde{Z}_{i}\ar[d]_{\widetilde{\pi}_{i}} & & \\
\widetilde{Z}_{i-1}\ar[d] & & \\
\cdot & & \\
\cdot\ar[d] & & \\
\widetilde{Z}_{1}\ar[d]_{\widetilde{\pi}_{1}} & & \\
\widetilde{Z}_{0} & & }
\end{equation*}
The next step in our proof refers to how to obtain the intersection form in $\widetilde{X}_{l-1}$ associated to the simple normal crossing divisor $\widetilde{D}_{\widetilde{X}_{l-1}}$. \\
By the projection formula (Proposition \ref{Pro3})
\begin{equation*}
\widetilde{H}_{\widetilde{X}_{l-1}, i_{1}}\cdot \widetilde{H}_{\widetilde{X}_{l-1}, i_{2}}\cdot\cdot\cdot \widetilde{H}_{\widetilde{X}_{l-1}, i_{d}}= \widetilde{f}_{l}^{*}(\widetilde{H}_{\widetilde{X}_{l-1}, i_{1}})\cdot f_{l}^{*}(\widetilde{H}_{\widetilde{X}_{l-1}, i_{2}})\cdot\cdot\cdot f_{l}^{*}(\widetilde{H}_{\widetilde{X}_{l-1}, i_{d}}),
\end{equation*}
Applying the result of Theorem \ref{ThmTotTrans} then
\begin{equation}\label{EqTotTrans}
\widetilde{H}_{\widetilde{X}_{l-1}, i_{1}}\cdot \widetilde{H}_{\widetilde{X}_{l-1}, i_{2}}\cdot\cdot\cdot \widetilde{H}_{\widetilde{X}_{s-1}, i_{d}}=(\widetilde{H}_{\widetilde{Z}_{s}, i_{1}}+\delta_{i_{1},i}\widetilde{H}_{\widetilde{Z}_{s},i})\cdot (\widetilde{H}_{\widetilde{Z}_{s}, i_{2}}+\delta_{i_{2},i}\widetilde{H}_{\widetilde{Z}_{s},i})\cdot\cdot\cdot (\widetilde{H}_{\widetilde{Z}_{s}, i_{d}}+\delta_{i_{d},i}\widetilde{H}_{\widetilde{Z}_{s},i}),
\end{equation}
where $\delta_{i_{j},i}=1$ if $\widetilde{H}_{\widetilde{Z}_{s},i}\cap \widetilde{H}_{\widetilde{Z}_{s},i_{j}}\neq\emptyset$ (see numerical characterization in Lemma \ref{Lem2}) and $\delta_{i_{j},i}=0$ otherwise.
\end{proof}

\setcounter{theorem}{12}

\begin{remark}\label{RemSRC}
It follows then that by iterating the above process, contracting a final divisor at each step, we will obtain a sequence of point blow-ups of length $l$. The algebraically marked sequence obtained depends on the choice of final components. Below we will prove that all the algebraically marked sequential morphisms associated to the sequences constructed in this way are algebraically equivalent. 
\end{remark}

\begin{proposition}\label{ProEqSM}
Any of the algebraically marked sequences obtained as in Remark \ref{RemSRC}, as composition of regular projective contractions from a fixed sky $Z$ and a fixed simple normal crossing divisor $E_{Z}$, are associated to algebraically marked sequential morphisms in the same algebraic equivalence class (see Definition \ref{DefAlEqSM}).
\end{proposition}

Before proving this, we need the following lemma.
\begin{lemma}
Fix the sky $Z$, fix a simple normal crossing divisor $E_{Z}$ and let us suppose that $H_{i}$ and $H_{j}$ are both finals for the sequential morphism $\pi: Z\rightarrow Z_{0}$. Then, there is an isomorphism $X_{m-2}\cong X_{m-2}^{'}$ making the following diagram commutative 
\begin{equation*}
\xymatrix{ & Z\ar[rd]^{f_{m}^{'}}\ar[ld]_{f_{m}} &  \\
X_{m-1}\ar[d]_{f_{m-1}} & & X_{m-1}^{'}\ar[d]^{f_{m-1}^{'}}  \\
X_{m-2}\ar[rr]^{\cong} & & X_{m-2}^{'}\ar[ll]}
\end{equation*}
where $f_{m}$ is the contraction of $H_{Z,i}$ and $f_{m-1}$ is the contraction of $H_{X_{m-1}, j}$, whereas $f_{m}^{'}$ is the contraction of $H_{Z, j}$ and $f_{m-1}^{'}$ is the contraction of $H_{X_{m-1}^{'}, i}$.
\end{lemma}

\begin{proof}
To begin with, if we denote by $O_{i, m-2}=f_{m-1}\circ f_{m}(H_{Z,i})$, $O_{j,m-2}=f_{m-1}(H_{X_{m-1}, j})$, $O_{j,m-2}^{'}=f_{m-1}^{'}\circ f_{m}^{'}(H_{Z,j})$ and $O_{i,m-2}^{'}=f_{m-1}^{'}(H_{X_{m-1}^{'}, i})$, then it follows that
\begin{equation*}
X_{m-2}\setminus\left\{O_{i, m-2}, O_{j, m-2}\right\}\cong Z\setminus H_{Z,i}\cup H_{Z,j}\cong X_{m-2}^{'}\setminus\left\{O_{i, m-2}^{'}, O_{j, m-2}^{'})\right\}
\end{equation*}
Let $W_{j}$ be an open affine open neighborhood of $O_{j, m-2}$. If we denote by $V_{j}$ to the inverse image $f_{m}^{-1}\circ f_{m-1}^{-1}(W_{j})$, then the image $W_{j}^{'}=f_{m}^{'}\circ f_{m-1}^{'}(V_{j})$ will be an affine open neighborhood of $O_{j, m-2}^{'}$. Then since $f_{m-1}\circ f_{m}\vert_{V_{j}}$ and $f_{m-1}^{'}\circ f_{m}^{'}\vert_{V_{j}}$ are both proper morphisms, it follows by Lemma \ref{Lem4} $W_{j}\cong W_{j}^{'}$. \\
If we denote by $W_{i}$ to an open affine neighborhood of $O_{i, m-2}$ and $W_{i}^{'}=f_{m-1}^{'}\circ f_{m}^{'}(V_{i})$, where $V_{i}$ is the inverse image $f_{m}^{-1}\circ f_{m-1}^{-1}(W_{i})$, then in a similar way we can prove that $W_{i}\cong W_{i}^{'}$, so it follows  $X_{m-2}\cong X_{m-2}^{'}$ since all isomorphisms are given by global sections.
\end{proof}

Consequently, we have the following corollary, which means that Proposition \ref{ProEqSM} holds for length $2$.
\begin{corollary}\label{EqSMl2}
If $Z$ is the sky of a sequence of point blow-ups of length $2$, then any of the two sequences of point blow-ups obtained following the procedure in \ref{RemSRC} are associated to algebraically marked sequential morphisms in the same algebraic equivalence class .
\end{corollary}

In order to prove Proposition \ref{ProEqSM} we need the following definition.
\begin{definition}
Consider two sequences of point blow-ups obtained as in Remark \ref{RemSRC}, obtained by composition of regular projective contractions from a fixed sky $Z$ and a fixed simple normal crossing divisor $E_{Z}$,
\begin{equation*}
\xymatrix{ 
Z\ar[r]^{f_{s}}\ar[d] & X_{s-1}\ar[r]^{f_{s-1}}\ar[d] & X_{s-2}\ar[r]^{f_{s-2}} & .\ar[r] & .\ar[r] & .\ar[r] & X_{2}\ar[r]^{f_{2}} & X_{1}\ar[r]^{f_{1}} & X_{0}\ar[d] \\
Z\ar[r]^{f_{s}^{'}}\ar[u] & X_{s-1}^{'}\ar[r]^{f_{s-1}^{'}}\ar[u] & X_{s-2}^{'}\ar[r]^{f_{s-2}^{'}} & .\ar[r] & .\ar[r] & .\ar[r] & X_{2}^{'}\ar[r]^{f_{2}^{'}} & X_{1}^{'}\ar[r]^{f_{1}^{'}} & X_{0}^{'}\ar[u]
}
\end{equation*}
We say that the sequences have the same end if at least the first contraction is common to both. i.e. one has $f_{s}=f^{'}_{s}$.
\end{definition}
\begin{proof}[Proof of Prop. \ref{ProEqSM}]
Let us suppose then that $Z$ is the sky of a sequence of $n+1$ point blow ups and that Proposition \ref{ProEqSM} is true for sequences of length lower or equal than $n$. If two sequences obtained as above $\rho:=f_{1}\circ f_{2}\circ...\circ f_{n}\circ f_{n+1}$ and $\rho^{'}:=f_{1}^{'}\circ f_{2}^{'}\circ...\circ f_{n}^{'}\circ f_{n+1}^{'}$ have the same end, then it is clear that both are associated to algebraically marked sequential morphism in the same algebraic equivalence class. It is a direct consequence of the fact that by hypothesis the assertion is true for sequences of length lower or equal than $n$.
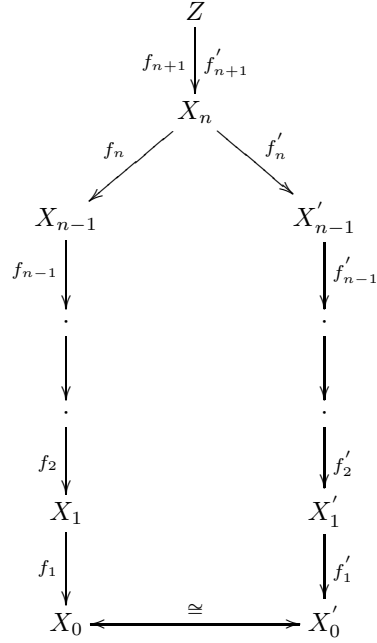
\begin{equation*}
\xymatrix{ & Z\ar[d]^{f_{n+1}^{'}}_{f_{n+1}} &  \\
 & X_{n}\ar[ld]_{f_{n}}\ar[rd]^{f_{n}^{'}} &   \\
X_{n-1}\ar[d]_{f_{n-1}}& & X_{n-1}^{'}\ar[d]^{f_{n-1}^{'}} \\
\cdot\ar[d] & & \cdot\ar[d] \\
\cdot\ar[d]_{f_{2}} & & \cdot\ar[d]^{f_{2}^{'}} \\
X_{1}\ar[d]_{f_{1}} & & X_{1}^{'}\ar[d]^{f_{1}^{'}} \\
X_{0}\ar[rr]^{\cong} & & X_{0}^{'}\ar[ll]}
\end{equation*}
If two sequences $\rho:=f_{1}\circ f_{2}\circ...\circ f_{n}\circ f_{n+1}$ and $\sigma:=g_{1}\circ g_{2}\circ...\circ g_{n}\circ g_{n+1}$ have not the same end, then let us suppose that $f_{n+1}$ and $g_{n+1}$ correspond to the contraction of $H_{Z,i}$ and $H_{Z,j}$ respectively. Consider all the sequences that belong to the tree contracting $H_{Z,j}$ first, there must exist a sequence $\rho^{'}:=f_{1}^{'}\circ f_{2}^{'}\circ...\circ f_{n}^{'}\circ f_{n+1}^{'}$ contracting $D_{X_{n},i}$ secondly. Analogously, if we consider all sequences contracting $H_{Z,i}$ first, there must exist a sequence $\sigma^{'}:=g_{1}^{'}\circ g_{2}^{'}\circ...\circ g_{n}^{'}\circ g_{n+1}^{'}$ contracting $D_{Y_{n},j}$ secondly.\\
By corollary \ref{EqSMl2} the sequences $f_{n}^{'}\circ f_{n+1}^{'}$ and $g_{n}^{'}\circ g_{n+1}^{'}$ of length $2$ are associated to sequential morphism in the same equivalence class, so it just remain to proof that $f_{1}^{'}\circ f_{2}^{'}\circ...\circ f_{n-2}^{'}\circ f_{n-1}^{'}$ belong to the same equivalence class that $g_{1}^{'}\circ g_{2}^{'}\circ...\circ g_{n-2}^{'}\circ g_{n-1}^{'}$. But this equivalence follows directly from the hypothesis, so $\rho^{'}\sim\sigma^{'}$.
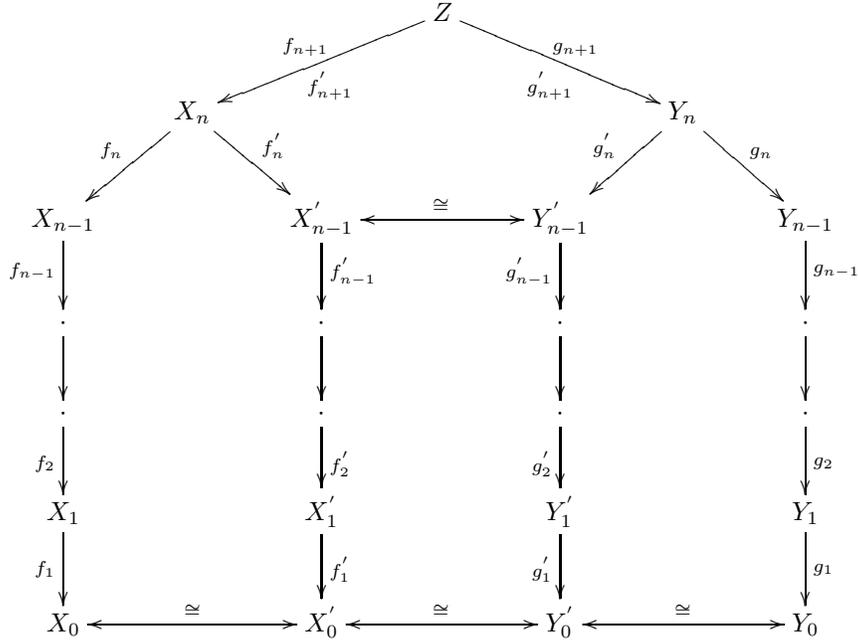
\begin{equation*}
\xymatrix{ & & & Z\ar[lld]_{f_{n+1}}^{f_{n+1}^{'}}\ar[rrd]^{g_{n+1}}_{g_{n+1}^{'}} & & &  \\
 & X_{n}\ar[ld]_{f_{n}}\ar[rd]^{f_{n}^{'}} & & & & Y_{n}\ar[ld]_{g_{n}^{'}}\ar[rd]^{g_{n}} & \\
X_{n-1}\ar[d]_{f_{n-1}} & & X_{n-1}^{'}\ar[d]^{f_{n-1}^{'}}\ar[rr]^{\cong} & & Y_{n-1}^{'}\ar[d]_{g_{n-1}^{'}}\ar[ll] & & Y_{n-1}\ar[d]^{g_{n-1}} \\
\cdot\ar[d] & & \cdot\ar[d] & & \cdot\ar[d] & & \cdot\ar[d] \\
\cdot\ar[d]_{f_{2}} & & \cdot\ar[d]^{f_{2}^{'}} & & \cdot\ar[d]_{g_{2}^{'}} & & \cdot\ar[d]^{g_{2}} \\
X_{1}\ar[d]_{f_{1}} & & X_{1}^{'}\ar[d]^{f_{1}^{'}} & & Y_{1}^{'}\ar[d]_{g_{1}^{'}} & & Y_{1}\ar[d]^{g_{1}} \\
X_{0}\ar[rr]^{\cong} & & X_{0}^{'}\ar[rr]^{\cong}\ar[ll] & & Y_{0}^{'}\ar[ll]\ar[rr]^{\cong} & & Y_{0}\ar[ll] }
\end{equation*}
Now since $\rho\sim\rho^{'}$ and $\sigma\sim\sigma^{'}$, then $\rho\sim\sigma$.
\end{proof}

\begin{proof}[Proof of Proposition \ref{Thm41}, Part $II$]
We can apply Proposition \ref{ProEqSM} to the morphism $\widetilde{Z}_{s}\rightarrow\widetilde{Z}_{0}$ an by scalar extension $\times_{Spec(\widetilde{K})} Spec(K)$ the algebraically marked sequences of point blow-ups constructed as above
\begin{equation*}
\xymatrix{
Z_{s}\ar[r]\ar[d] & X_{l-1}\ar[r]\ar[d] & \cdot\ar[r] & X_{1}\ar[r]\ar[d] & X_{0}\ar[d] \\
\widetilde{Z}_{s}\ar[r] & \widetilde{X}_{l-1}\ar[r] & \cdot\ar[r] & \widetilde{X}_{1}\ar[r] & \widetilde{X}_{0}
}
\end{equation*}
where $X_{i}\cong \widetilde{X}_{i}\times_{Spec(\widetilde{K})} Spec(K)$, so the theorem is proved.
\end{proof}

Now we are ready to prove our main theorem.
\begin{theorem}\label{Thm21}
Let $(\pi: Z_{s}\rightarrow Z_{0}, \sqcup_{i=1}^{l}F_{i})_{alg}$ and $(\pi^{'}: Z_{s}^{'}\rightarrow Z_{0}^{'}, \sqcup_{i=1}^{l^{'}}F^{'}_{i})_{alg}$ be two algebraically marked sequential morphisms. Then  they are algebraically equivalent over $K$ as in Definition \ref{DefAlEqSM} if and only if
 there exist algebraically marked sequences of point blow-ups $(Z_{s},...,Z_{0},\sqcup_{i=1}^{l}F_{i},\pi)_{alg}$ and $(Z^{'}_{s},...,Z^{'}_{0},\sqcup_{i=1}^{l^{'}}F_{i}^{'},\pi^{'})_{alg}$  associated to $(\pi: Z_{s}\rightarrow Z_{0}, \sqcup_{i=1}^{l}F_{i})_{alg}$ and $(\pi^{'}: Z_{s}^{'}\rightarrow Z_{0}^{'}, \sqcup_{i=1}^{l^{'}}F^{'}_{i})_{alg}$ respectively such that they are algebraically equivalent over $K$ as in Definition \ref{DefAlEqS}.
\end{theorem}

\begin{proof}
If two algebraically marked sequences of point blow-ups are algebraically equivalent, then it follows directly by Definition \ref{DefAlEqS} that the associated algebraically marked sequential morphisms are algebraically equivalent too. \\
Now we will prove that if two algebraically marked sequential morphism $(\pi: Z_{s}\rightarrow Z_{0}, \sqcup_{i=1}^{l}F_{i})_{alg}$ and $(\pi^{'}: Z_{s}^{'}\rightarrow Z_{0}^{'}, \sqcup_{i=1}^{l^{'}}F^{'}_{i})_{alg}$ are algebraically equivalent, then there exist algebraically marked sequences of point blow-ups associated to them that are algebraically equivalent too. By Proposition \ref{Thm41} given a certain sky $Z$ associated to an algebraically marked sequential morphism $(\pi: Z_{s}\rightarrow Z_{0}, \sqcup_{i=1}^{l}F_{i})_{alg}$ , all the algebraically marked sequences of point blow-ups obtained by regular projective contractions are associated to algebraically marked sequential morphisms in the same algebraic equivalence class. Since $(\pi: Z_{s}\rightarrow Z_{0}, \sqcup_{i=1}^{l}F_{i})_{alg}$ and $(\pi^{'}: Z_{s}^{'}\rightarrow Z_{0}^{'}, \sqcup_{i=1}^{l^{'}}F^{'}_{i})_{alg}$ are algebraically equivalent, then $\exists\widetilde{K}\subset K$ such that there exists an isomorphism $\widetilde{b}: \widetilde{Z}\rightarrow \widetilde{Z^{'}}$. By applying Proposition \ref{Pro8} and Proposition \ref{ProEqSM} we conclude the result by scalar extension $\times_{Spec(\widetilde{K})} Spec(K)$.  
\end{proof}

\subsection{Examples and comments}

In the following example we will show that the numerical information that appears in the combinatorial object associated to a sequential morphism, that is, the multilinear map $\Phi_{Z, \left(\sqcup_{i=1}^{l}F_{i}\right)_{Z, \pi}, K}$, is minimal in some sense for our classification purposes (it is not possible to establish the bijection $\Theta: (\Phi_{Z, \left(\sqcup_{i=1}^{l}F_{i}\right)_{Z, \pi}, K})/\sim\rightarrow (T_{(Z_{s},...,Z_{0}, \sqcup_{i=1}^{l}F_{i}, \pi)})/\sim$ otherwise)

\begin{example}\label{Exmp1}
Consider a $d-$dimensional smooth algebraic variety $X$ over $k$, where $k$ is a perfect field not algebraically closed, and $d$ is an even natural number. Let $P=P_{1}$ be a $k-$rational point on $X$ and let $\pi_{1}:Z_{1}\rightarrow X$ be the blow-up with center $P$. Take a point $P_{2}$ on $Z_{1}$ such that the degree of the field extension $k\subset k(P_{2})$ is $1$, and a point $Q_{2}$ in $Z_{1}$ such that $\left[k(Q_{2}):k\right]=2$. Now, let $\pi_{2}:Z_{2}\rightarrow Z_{2}$ (respectively $\pi_{2}^{'}: Z_{2}^{'}\rightarrow Z_{1}$) be the blowing up of $P_{2}$ (of $Q_{2}$, respectively) and take a point $P_{3}$ in $Z_{2}$ such that $P_{3}$ is not proximate to $P_{1}$ and such that the degree of the field extension $k(P_{2})\subset k(P_{3})$ is $3$ (respectively, a point $Q_{3}$ in $Z_{2}^{'}$ no proximate to $P_{1}$ such that $\left[k(Q_{3}):k(Q_{2})\right]=1$). Consider now the blow-up $\pi_{3}: Z_{3}\rightarrow Z_{2}$ ($\pi_{3}^{'}: Z_{3}^{'}\rightarrow Z_{2}^{'}$) of $Z_{2}$ ($Z_{2}^{'}$) with center $P_{3}$ ($Q_{3}$ respectively).
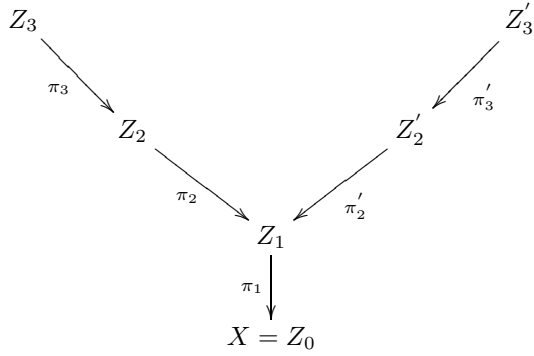
\begin{equation*}
\xymatrix{ Z_{3}\ar[rd]_{\pi_{3}} & & & & Z_{3}^{'}\ar[ld]^{\pi_{3}^{'}} \\
 & Z_{2}\ar[rd]_{\pi_{2}} & & Z_{2}^{'}\ar[ld]^{\pi_{2}^{'}} & \\
 & & Z_{1}\ar[d]_{\pi_{1}} & & \\
 & & X=Z_{0} & &}
\end{equation*}
Let $\Delta: E\rightarrow\overbrace{E\times E\times\cdots\times E}^{d}$ be the diagonal map, that is $\Delta\left(E_{i}\right)=\left(E_{i}, E_{i}, \ldots, E_{i}\right)$. Instead of the multilinear map $\Phi_{Z_{3}, \left(\sqcup_{i=1}^{3}F_{i}\right)_{Z_{3}, \pi}, k}$ ($\Phi_{Z_{3}^{'}, \left(\sqcup_{i=1}^{3}F_{i}^{'}\right)_{Z_{3}^{'}, \pi^{'}}, k})$, let us consider the multilinear map associated to $\mathcal{I}_{Z_{3}, \left(\sqcup_{i=1}^{3}F_{i}\right)_{Z_{3}, \pi}, k}\circ\Delta$ ($\mathcal{I}_{Z_{3}^{'}, \left(\sqcup_{i=1}^{3}F_{i}^{'}\right)_{Z_{3}^{'}, \pi^{'}}, k}\circ\Delta$ respectively) we will refer to as the diagonal of the multilinear map. Then, the numerical information contained in this combinatorial object for the sequential morphisms $\pi=\pi_{1}\circ\pi_{2}\circ\pi_{3}$ and $\pi^{'}=\pi_{1}^{'}\circ\pi_{2}^{'}\circ\pi_{3}^{'}$ is the following one, respectively:
\begin{equation*}
\left\{
\begin{aligned}
(1,0,0)^{d}&=-2, \\
(0,1,0)^{d}&= -4, \\
(0,0,1)^{d}&=-3;
\end{aligned}
\right.
\end{equation*}
\begin{equation*}
\left\{
\begin{aligned}
(1,0,0)^{d}&=-3 \\
(0,1,0)^{d}&=-4 \\ 
(0,0,1)^{d}&=-2 
\end{aligned}
\right.
\end{equation*}
Now, if we had just considered the multilinear map associated to diagonal of the $d-$ary intersection form, then by naturally extending our equivalence notion (see Definition \ref{ComEqSM}),  it follows that these multilinear maps would be equivalent. However, it is clear that the weighted directed graphs associated to these sequences of point blow-ups, $\left(Z_{3}, \ldots, Z_{0}, \pi\right)$ and $\left(Z_{3}^{'}, \ldots, Z_{0}, \pi^{'}\right)$ are not equivalent:
\begin{align*}
T_{(Z_{3},Z_{2},Z_{1},Z_{0}, \sqcup_{i=1}^{3}F_{i}, \pi)}\equiv & \begin{cases}\xymatrix{ 3\ar[d] \\ 1\ar[d] \\ 1},\end{cases} \\
T_{(Z_{3}^{'},Z_{2}^{'},Z_{1},Z_{0}, \sqcup_{i=1}^{3}F_{i}^{'}, \pi^{'})}\equiv & \begin{cases}\xymatrix{ 2\ar[d] \\ 2\ar[d] \\ 1}. \end{cases}
\end{align*}
This is due to the fact that although the multilinear maps associated to the diagonals are equivalent, the whole multilinear maps are not:
\begin{equation*}
\left\{
\begin{aligned}
(1,0,0)^{r}\cdot (0,1,0)^{s}&=(-1)^{r+1}, \\
(0,1,0)^{r}\cdot (0,0,1)^{s}&=(-1)^{r}\cdot(-3),
\end{aligned}
\right.
\end{equation*}
for the multilinear map $\Phi_{Z_{3}, \left(\sqcup_{i=1}^{3}F_{i}\right)_{Z_{3}, \pi}, k}$, and
\begin{equation*}
\left\{
\begin{aligned}
(1,0,0)^{r}\cdot (0,1,0)^{s}&=(-1)^{r+1}\cdot 2, \\
(0,1,0)^{r}\cdot (0,0,1)^{s}&=(-1)^{r}\cdot(-2),
\end{aligned}
\right.
\end{equation*}
for the multilinear map $\Phi_{Z_{3}^{'}, \left(\sqcup_{i=1}^{3}F_{i}^{'}\right)_{Z_{3}^{'}, \pi^{'}}, k}$.
\end{example}

\begin{remark}
The above result makes a huge difference with respect to the particular case where we just consider sequences of point blow-ups defined over an algebraically closed field. In that case, if the multilinear maps associated to the diagonals are equivalent then the whole multilinear maps are also equivalent.
\end{remark}

\subsection*{Acknowledgments}

We would like to thank Professor A. Campillo for stimulating discussions and his valuable comments and suggestions on the article.

\end{document}